\documentclass[11pt]{article}

\usepackage[english]{babel}

\usepackage[letterpaper,top=2cm,bottom=2cm,left=2cm,right=2cm,marginparwidth=1.75cm]{geometry}

\usepackage{amsmath}
\usepackage{mathtools}
\mathtoolsset{showonlyrefs=false}
\usepackage{amssymb}
\usepackage{graphicx}
\usepackage{tikz}
\usepackage{tikz-3dplot}
\usepackage{pgfplots}
\pgfplotsset{compat=1.16}
\usepackage{amsmath,amssymb,graphicx,algorithm,algpseudocode}
\usepackage{mathrsfs}
\usepackage{ulem}
\mathtoolsset{showonlyrefs=true}
\usetikzlibrary{positioning}

\usepackage[algo2e]{algorithm2e}

\usepackage[colorlinks=true, allcolors=blue]{hyperref}
\usepackage{color}
\usepackage{here}
\usepackage{epsfig}
\usepackage{ulem}
\usepackage[justification=centering]{caption}

\usepackage{amsthm}
\theoremstyle{plain}
\newtheorem{theorem}{Theorem}[section]
\theoremstyle{definition}

\newtheorem{lemma}[theorem]{Lemma}

\newtheorem{remark}{Remark}[section]
\newtheorem{conj}[theorem]{Conjecture}

\newtheorem*{theorem4.1}{Theorem 4.1}
\newtheorem*{theorem1.2}{Theorem 1.2 of \cite{endo2025second1}}

\newcommand{\CG}{{\mbox{\tiny CG}}}

\usepackage{color}

\newcommand{\ENDO}[1]{{\color{blue}{#1}}}
\renewcommand{\ENDO}[1]{#1}

\title{The Second Dirichlet Eigenvalue is Simple on Every Non-equilateral Triangle, Part II: Nearly Equilateral Triangles}
\author{Ryoki Endo\thanks{Graduate School of Science and Technology, Niigata University, Niigata, Japan
(endo@m.sc.niigata-u.ac.jp).}, ~~Xuefeng Liu\thanks{Department of Information and Sciences, Tokyo Woman's Christian University, Tokyo, Japan (xfliu@cis.twcu.ac.jp) (Corresponding author).}}

\begin{document}

\maketitle

\abstract{ 
This paper solves the open problem of the simplicity of the second Dirichlet eigenvalue for nearly equilateral triangles, offering a complete solution to Conjecture 6.47 posed by R. Laugesen and B. Siudeja in A. Henrot's book ``Shape Optimization and Spectral Theory."
Our proof is achieved by introducing a new difference quotient formula for the behavior of nearly degenerate eigenvalues resulting from domain perturbations, and a novel numerical algorithm that rigorously estimates it using verified computation.
}

\section{Introduction}

\begin{figure}[H] 
        \centering
        \includegraphics[keepaspectratio, scale=0.40]{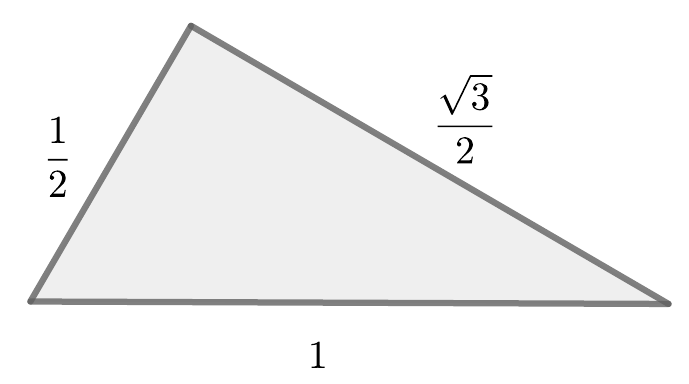}
        \caption{Half equilateral triangle}
        \label{fig:half-equilateral}
    \end{figure}

Determining the eigenvalue multiplicity of the Laplace operator over a given domain is generally challenging. By evaluating the lower and upper bounds of eigenvalues, one can validate the simplicity of the eigenvalues if they are separated from each other. In \cite{liu2013verified,liu2015framework}, the second author developed an algorithm to provide guaranteed lower and upper bounds for the Laplacian eigenvalues over polygonal domains of general shapes, where the finite element method (FEM) plays an important role. These methods can determine the simplicity of eigenvalues in case they are well separated from the rest. For example, the second Dirichlet eigenvalue is simple over the concrete triangular domain as shown in Figure \ref{fig:half-equilateral}, since the following eigenvalue bounds are obtained by the method proposed in \cite{liu2015framework}:
    \begin{equation}
        122.63\leq \lambda_1\leq 123.04,~227.43\leq \lambda_2\leq  228.79,~332.03\leq \lambda_3\leq 334.81
    \end{equation}
    However, if eigenvalues are nearly degenerate, it becomes difficult to separate them using eigenvalue bounds. For example, it is known that $\lambda_2=\lambda_3$ for the Dirichlet eigenvalues over equilateral triangles; see \cite{mccartin2003eigenstructure}. Due to the continuity of eigenvalues with respect to domain perturbation, $\lambda_2$ and $\lambda_3$ form a tight cluster over nearly equilateral triangles, making it extremely challenging to determine the multiplicities of these two eigenvalues. The following open conjecture \ENDO{posed by R. Laugesen and B. Siudeja}, which will be discussed in this paper, addresses the multiplicity of Dirichlet eigenvalues for triangular domains.

    \medskip
    \begin{conj}[Conjecture 6.47 of \cite{henrot2017shape}]\label{len:main-conjecture}
	The second Dirichlet eigenvalue is simple for every non-equilateral triangle.
    \end{conj}
    \medskip

    Despite the substantial contributions from the vast mathematical community, numerous questions regarding the eigenvalues of triangles remain unresolved. 
    Except for Conjecture 1, below, we present several conjectures and open problems related to triangular domains. For further details, see \cite{henrot2017shape}. \begin{enumerate} \item Can the shape of a triangular drum be discerned from its first three tones? That is, can the geometry of a triangle be uniquely determined by its first three Dirichlet eigenvalues? \item The equilateral triangle is conjectured to minimize the $k$-th Dirichlet eigenvalue ($k \geq 4$) among all triangles with a fixed diameter. 
    %\item The second Dirichlet eigenvalue is expected to be simple for every non-equilateral triangle. 
    \item Among $N$-gons ($N = 3, 4, \dots$) with a given area, the regular $N$-gon is conjectured to maximize the second Neumann eigenvalue. \end{enumerate}

    In our companion paper \cite{endo2025second1}, we provided a partial result confirming Conjecture \ref{len:main-conjecture} for the case of collapsing triangles:
    \begin{theorem1.2}
    The second Dirichlet eigenvalue is simple for every triangle with its minimum normalized height\footnote{The minimum normalized height of a triangle is the height measured relative to its longest side, with the triangle scaled such that the longest side has unit length.}  less than or equal to 
    $\tan(\pi/60)/2$.
    \end{theorem1.2}
%     \begin{remark}
%         For a narrow triangle, the eigenvalues diverge to infinity as the height of the triangle approaches zero. 
% In \cite{ourmieres2015dirichlet}, Ourmières--Bonafos established the following asymptotic expansion for the Dirichlet eigenvalues over the triangle \( T(s,t) \) with vertices $(-1,0),(1,0),(s,t)$ :
% \begin{equation}\label{eq:bonafos}
% \lambda_k^t \sim t^{-2} \left( \pi^2 + (2\pi^2)^{2/3} \kappa_k t^{2/3} + \beta_{3,k} t + \cdots \right).    
% \end{equation}
% Here, \( \kappa_k \) $(k=1,2,\cdots)$ are the eigenvalues of the Schrödinger operator
% \[l_s^{\text{mod}}:=
% -\frac{d^2}{dx^2} + \left(\frac{1}{1+s}\mathbf{1}_{\mathbb{R}_-}+\frac{1}{1-s}\mathbf{1}_{\mathbb{R}_+}\right)|x|
% \]
% defined on \( H^1(\mathbb{R}) \).
% \end{remark}
    This paper proves the conjecture for the complementary case, establishing the following result:
    \ENDO{\begin{theorem4.1}
    The second Dirichlet eigenvalue is simple for every non-equilateral triangle with its minimum normalized height \footnote{The minimum normalized height of a triangle is the height measured relative to its longest side, with the triangle scaled such that the longest side has unit length.} greater than or equal to $\tan(\pi/60)/2$.
    \end{theorem4.1}}
    
    \medskip
    
    Combining the result above with Theorem~1.2 from our companion paper \cite{endo2025second1}, we obtain a complete proof of Conjecture \ref{len:main-conjecture}:
    \begin{theorem}\label{thm:main-conjecture}
The second Dirichlet eigenvalue is simple for every non-equilateral triangle.
\end{theorem}
    
\medskip        
    
    In this paper, we propose a method to determine the multiplicities of nearly degenerate eigenvalues. For this purpose, we study the difference quotient with respect to domain perturbation as defined in \eqref{def:difference-quotient}, and we derive the difference quotient formula in Theorem \ref{lem:Fte-basis-eigen}. The rigorous estimation of the difference quotient in Theorem \ref{lem:M-star-N-star-estimation} is obtained by further applying the error estimation of eigenspaces established in \cite{liu2022fully}.
\\  
    As an application of Theorem \ref{lem:Fte-basis-eigen} and Theorem \ref{lem:M-star-N-star-estimation}, we investigate the multiplicity of Dirichlet eigenvalues over nearly equilateral triangles. As a result, we partially solve Conjecture \ENDO{\ref{len:main-conjecture}} and derive the following result as shown in Theorem \ref{len:partial-result}:
        \medskip

\medskip

    To obtain the local behavior of eigenvalues with respect to domain perturbations, it is effective to analyze the first-order variation of eigenvalues. Fundamental results on the first-order variation of eigenvalues include Hadamard's shape derivative formula in \cite{hadamard1908memoire} and Rousselet's results in \cite{rousselet1983shape} on the directional derivative of multiple eigenvalues from the 1980s. The analysis on derivatives of eigenvalues reveals the local behavior of eigenvalues with respect to domain perturbations. For example, in \cite{endo2023shape}, we investigated the local behavior of eigenvalues over nearly equilateral triangles by rigorously evaluating the shape derivative of the first Dirichlet eigenvalue, and used this information to solve a certain shape optimization problem over triangular domains. However, when eigenvalues are nearly degenerate and their multiplicities are unknown, the analysis on derivatives does not directly provide information about the behavior of eigenvalues over a neighborhood with a given radius. In our opinion, the existing theorems in \cite{hadamard1908memoire} and \cite{rousselet1983shape} on eigenvalue derivatives are not enough to study the multiplicities involved in Conjecture \ref{len:main-conjecture}.
\\

To better understand the limitations of using derivatives in the case of nearly degenerate eigenvalues and the motivation behind our approach, let us recall Rousselet's theorem on the directional derivative of multiple eigenvalues \cite{rousselet1983shape}.

\medskip

For an open set $\Omega_0$ of $\mathbb{R}^2$, consider a family of functions $\Phi(t)$ that satisfies the following conditions:
\begin{equation*}
	\Phi(t)\in W^{1,\infty}(\mathbb{R}^2,\mathbb{R}^2) \mbox{ is differentiable at $0$ with }\Phi(0)=I,\Phi'(0)=V,
		\end{equation*}
		where $W^{1,\infty}(\mathbb{R}^2,\mathbb{R}^2)$ is the set of bounded Lipschitz maps from $\mathbb{R}^2$ to itself, $I$ is the identity, and $V$ is a vector field. For the simplicity of notation, let $\Omega_t=\Phi(t)(\Omega_0)$ for $t\in[0,T)$. Let $\lambda_i(\Omega_t)$ be the $i$-th Dirichlet eigenvalue over $\Omega_t$.
  
  To study the variation of tightly clustered eigenvalues, let us consider the difference quotient of eigenvalues defined by
        \begin{equation}\label{def:difference-quotient}
            D_t\lambda_i:=\frac{\lambda_i(\Omega_t)-\lambda(\Omega_0)}{t}.
        \end{equation}	

        Let us quote the result of Hadamard  \cite{hadamard1908memoire} and Rousselet  \cite{rousselet1983shape} on the shape derivative of eigenvalues.
			\begin{theorem}{\cite[Theorem 3.2]{rousselet1983shape}}\label{lem:rousselet}
				Let $\lambda(\Omega_0)$ be a Dirichlet eigenvalue of the Laplacian over $\Omega_0$ with multiplicity $m\in\mathbb{N}$. Suppose that, for $t$ small enough, there are exactly $m$ eigenvalues $\lambda_1(\Omega_t),\lambda_2(\Omega_t)\cdots,\lambda_m(\Omega_t)$ in the neighborhood of $\lambda(\Omega_0)$ in $\mathbb{R}$.
				Let $\tilde a',\tilde b'$ be the Fréchet derivative of the bilinear forms $(\nabla\cdot,\nabla\cdot)_{L^2}$ and $(\cdot,\cdot)_{L^2}$ that define the weak form of Dirichlet eigenvalue problems. For the concrete definitions, see equations (32) and (107) in \cite{rousselet1983shape}.
\\												
				If $m=1$ (simple eigenvalue), then $\lambda(\Omega_t)$ is differentiable at $t=0$ and 
				\begin{equation*}
					\lambda'(\Omega_t)=\tilde a'(u,u)-\lambda(\Omega_0)\tilde b'(u,u),
				\end{equation*}
				where $u$ is an $L^2$-normalized eigenfunction associated with $\lambda(\Omega_0)$.
\\												
				If $m>1$ (multiple eigenvalue), then $\lambda_{i}(\Omega_t)$ is directionally differentiable. Let
				\begin{equation*}
					\lambda_{i,+}'(\Omega_t):=\lim_{\substack{t\to 0+}}D_t\lambda_i.
				\end{equation*}
				Moreover, directional derivatives $\lambda_{i,+}'(\Omega_t)$ $(i=1,\cdots,m)$ are the eigenvalues of the matrix $M$ of entries
				\begin{equation*}
					M_{ij}=\tilde a'(u_i,u_j)-\lambda(\Omega_0)\tilde b'(u_i,u_j),~~i,j=1,\cdots,m,
				\end{equation*}
				where $u_1,\cdots,u_m$ are $L^2$-orthonormal eigenfunctions associated with $\lambda(\Omega_0)$.
			\end{theorem}
%                 Theoretically, the second and third eigenvalues over a nearly equilateral triangle, i.e., $\Omega_t$ ($t\in(0,\delta]$) for small enough $\delta>0$, can be separated by evaluating the directional derivative of the eigenvalues over the equilateral triangle. 
%                 However, for practical problems, one cannot determine the concrete value of such a $\delta$
% \footnote{Let $\Omega_0$ be the equilateral triangle with the unit radius. Note that $\lambda_2(\Omega_0)=\lambda_3(\Omega_0)$ due to the equilateral shape of $\Omega_0$.
% Let $\alpha:=\lim_{t\to 0+}D_t\lambda_2$ and $\beta:=\lim_{t\to 0+}D_t\lambda_3$. 
% In case of $\alpha<\beta$, theoretically there exists $\delta>0$ such that
% $D_t\lambda_2< (\alpha+\beta)/2 < D_t\lambda_3$ for all $t\in (0,\delta].$
% However, this argument does not provide an explicit value of $\delta$.}.

\medskip

\ENDO{The above formula for simple or multiple eigenvalues is not applicable when the multiplicity itself is unknown. Furthermore, even for eigenvalues known a priori to be simple, if they are clustered, the numerical computation of the corresponding eigenfunctions is an ill-posed problem. This instability makes it prohibitive to use Hadamard's formula, which requires the true eigenfunctions, to accurately compute shape derivatives.

To overcome these obstacles, we introduce a difference quotient formula that provides a robust method for analyzing eigenvalue behavior without relying on the direct computation of unstable eigenfunctions. This formula, detailed in Theorem \ref{lem:Fte-basis-eigen}, allows for the rigorous separation of the difference quotients, providing a new strategy to prove eigenvalue simplicity in these challenging scenarios. Using this formula, we can separate the eigenvalues by validating the following relation:}

\begin{equation}\label{eq:epsilon-delta}
D_t\lambda_2 \in [\underline{\alpha},\overline{\alpha}],~~
D_t\lambda_3 \in [\underline{\beta},\overline{\beta}]~~\forall t\in (0,\delta];~\overline{\alpha}<\underline{\beta}.
\end{equation}

\ENDO{
Let us summarize the features of the methods proposed in this paper:

\begin{enumerate}
    \item[(a)] Compared with existing theories on shape derivatives relying on a priori information about the simplicity, our proposed method can deal with cases where the simplicity is unknown when a perturbation is applied to the domains with multiple eigenvalues.
    \item[(b)]This formula enables the stable computation of the eigenfunctions associated with tightly clustered eigenvalues that are difficult to obtain directly. This application is detailed in our recent preprint \cite{endo2025stable}, where we numerically demonstrate that the method recovers eigenfunctions with their correct theoretical symmetries, even under extremely small domain perturbations.
\end{enumerate}

}

\medskip

	The remainder of this paper is organized as follows. In Section 2, we review the error estimation of eigenfunctions, which will be utilized in the estimation of the difference quotient formula for eigenvalues. Section 3 presents the derivation of the difference quotient formula and its estimation. In Section 4, we provide a computer-assisted proof of Theorem \ref{len:partial-result}. The code for the computer-assisted proof is publicly available at the GitHub repository \url{https://github.com/ryendo/DirichletSimplicityClustered}.

	\section{Preliminary}\label{section:preliminary}
	We start our discussion using the standard notation of Sobolev spaces. Let $T \subset \mathbb{R}^2$ be a triangular domain. Denote by $L^2(T)$ the space of all real-valued square-integrable functions on $T$, and by $H^1(T)$ the Sobolev space of functions in $L^2(T)$ whose weak derivatives are also in $L^2(T)$. Furthermore, let $H^1_0(T)$ be the subspace of $H^1(T)$ consisting of functions that vanish on the boundary of $T$. The $L^2$-norm of $v \in L^2(T)$ is denoted by $\|v\|_{T}$, and the inner product in $L^2(T)$ or $(L^2(T))^2$ is represented by $(\cdot,\cdot)_T$. Let $\nabla$ denotes the gradient operator for functions in $H^1(T)$. It is worth noting that $(\nabla\cdot,\nabla\cdot)_T$ defines an inner product on $H^1_0(T)$ due to the imposed boundary conditions.

The weak formulation of the Dirichlet Laplacian eigenvalue problem $-\Delta u = \lambda u$ is to find $u \in H^1_0(T) \setminus \{0\}$ and  $\lambda > 0$ such that 
\begin{equation}
\label{eq:eigenvalue-problem}
 (\nabla u, \nabla v)_T = \lambda (u, v)_T \quad \forall v \in H^1_0(T).
\end{equation}
 Since the inverse of the Laplacian is a compact and self-adjoint operator, the spectral theorem implies that the problem \eqref{eq:eigenvalue-problem} has a countably infinite spectrum of eigenvalues, which can be ordered as $0 < \lambda_1(T) < \lambda_2(T) \leq \lambda_3(T) \leq \cdots$.

	\medskip
								
	Let us introduce distances between subspaces of $H^1_0(T)$ to evaluate the error of numerical approximation of eigenspaces, as well as the perturbation of the eigenspace due to domain perturbations.  Given two subspaces $E$ and $E^h$ of $H^1_0(T)$, the directed distances 
	$\delta_a,\delta_b$, $\bar{\delta}_a$ and  $\bar{\delta}_b$  are defined by
	\begin{gather*}
		\delta_a(E,E^h):=\max_{\substack{v\in E\\\|\nabla v\|_T=1}}\min_{v_h\in E^h}\|\nabla v-\nabla v_h\|_T,~~ 
		\delta_b(E,E^h):=\max_{\substack{v\in E\\\|v\|_T=1}}\min_{v_h\in E^h}\|v-v_h\|_T,\\
		\bar{\delta}_a(E,E^h):=\max_{\substack{v\in E\\\|v\|_T=1}}\min_{\substack{v_h\in E^h\\\|v_h\|_T=1}}\|\nabla v-\nabla v_h\|_T,~~
		\bar{\delta}_b(E,E^h):=\max_{\substack{v\in E\\\|v\|_T=1}}\min_{\substack{v_h\in E^h\\\|v_h\|_T=1}}\|v-v_h\|_T.
	\end{gather*}
								
\begin{figure}[H]
		\centering
		\includegraphics[keepaspectratio, scale=0.18]{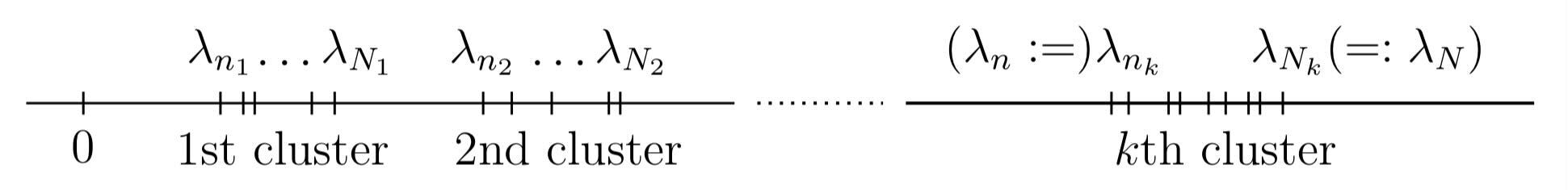}
		\caption{\label{fig:eigenvalue-distribution} Clusters of eigenvalues}
	\end{figure}

	In order to formulate the error of eigenfunctions, let us introduce a notation for clusters of eigenvalues. Let $n_k$ and $N_k$
	denote the indices of the first and the last eigenvalues in the $k$th cluster; see
	Figure \ref{fig:eigenvalue-distribution}. Note that eigenvalues in a cluster need not be equal to each other.
	We consider the $k$th cluster of interest, and set $n = n_k$ and $N = N_k$ to
	simplify the notation. 
								
	Let $E_k$ be the space of exact eigenfunctions associated to $k$th cluster of eigenvalues: 
	\begin{equation}\label{eq:def-Ek}
		E_k=\mbox{span}\{u_{n},u_{n+1},\cdots,u_{N}\}.
	\end{equation}
	Similarly, approximations $ u_{n,h} \in H^1_0(T)$ of exact eigenfunctions $u_i$ ($i=n,\cdots,N$) form the corresponding approximate space: 
	\begin{equation}\label{eq:def-Ehatk}
		E^h_{k}=\mbox{span}\{u_{n,h},u_{n+1,h},\cdots,\ u_{N,h}\}.
	\end{equation}
	Let $\lambda_{n,h}:=\min_{u_h\in E^h_k}\|\nabla u_h\|^2_T/\| u_h\|_T^2$ and $\lambda_{N,h}:=\max_{u_h\in E^h_k}\|\nabla u_h\|^2_T/\|u_h\|_T^2$.
	% Let us introduce measures of non-orthogonality between finite dimensional subspaces $E$ and $E'$ of $V(\mathcal{C})$ both in the $a$- and $b$-sense as
	% \begin{equation*}
	%     \hat\varepsilon_a(E,E')=\max_{\substack{v\in E\\\|\nabla v\|=1}}\max_{\substack{v'\in E'\\ \|\nabla v'\|=1}}(\nabla v,\nabla v'),~~\mbox{and} ~~\hat\varepsilon_b(E,E')=\max_{\substack{v\in E\\\|v\|=1}}\max_{\substack{\gamma v'\in E'\\ \|v'\|=1}}(v,v').
	% \end{equation*}
	Let us introduce measures of non-orthogonality between $E$ and $E^h$ by 
	\begin{equation*}
		\varepsilon^h_a(E,E^h)=\max_{\substack{v\in E\\\|\nabla v\|_T=1}}\max_{\substack{v_h\in E^h\\ \|\nabla v_h\|_T=1}}(\nabla v,\nabla v_h)_T,~~\varepsilon^h_b(E,E^h)=\max_{\substack{v\in E\\\|v\|_T=1}}\max_{\substack{v_h\in E^h\\ \|v_h\|_T=1}}(v,v_h)_T.
	\end{equation*}

 Let us quote the following well known formula from [\cite{boffi2010finite}, page 55] for the exact eigenpair $(\lambda_i,u_i)$ and its approximation $(\lambda_{i,h},u_{i,h})$:
        \begin{equation}\label{eq:fundamental-formula}
            \left\|\nabla u_i-\nabla u_{i,h}\right\|_T^2=\lambda_i\left\|u_i-u_{i,h}\right\|_T^2-\left(\lambda_i-\lambda_{i,h}\right)\left\|u_{i,h}\right\|_T^2~(i=n,\cdots,N).
        \end{equation}
 
	The following estimation will be used to measure the error between the clusters of eigenfunctions. As an important feature of this estimation, it gives error estimate of eigenspace approximation by only using the information of eigenvalues and the orthogonality of approximate eigenfunctions. \ENDO{Note that all the eigenvalues discussed in this paper can be estimated either directly using the recently developed estimation method (see Lemmas \ref{lem:est-tau} and \ref{thm:Lehmann-Goerisch}) or through perturbation estimates (see Lemma \ref{lem:eig-perturbation}).}

	\begin{lemma}{[Theorem 1 of Liu-Vejchodsk{\`y} \cite{liu2022fully}]}\label{lem:eigenvec-estimation-algorithm-1}
		Let $\rho>0$ be a quantity such that $\lambda_{n_k}<\rho\leq\lambda_{N_k+1}$. Then, we have
		\begin{align}\label{eq:estimation-for-delta}
		\begin{split}
            \delta^2_a(E_k,E^h_{k})&\leq\frac{\rho(\lambda_{N_k,h}-\lambda_{n_k})+\lambda_{n_k}\lambda_{N_k,h}\theta_a^{(k)}}{\lambda_{N_k,h}(\rho-\lambda_{n_k})}(=:\mathcal{E}^2_a(E_k,E^h_k)),\\
            \delta^2_b(E_k,E^h_{k})&\leq\frac{\lambda_{N_k,h}-\lambda_{n_k}+\theta_b^{(k)}}{\rho-\lambda_{n_k}}(=:\mathcal{E}^2_b(E_k,E^h_k)).
            \end{split}
		\end{align}
		Here, 
		\begin{align}
            \begin{split}
			\theta_a^{(k)}&=\sum_{l=1}^{k-1}\frac{\rho-\lambda_{n_l}}{\lambda_{n_l}}\left[\varepsilon^h_a(E^h_{l},E^h_{k})+\delta_a(E_{l},E^h_{l})\right]^2,\\
            \theta_b^{(k)}&=\sum_{l=1}^{k-1}(\rho-\lambda_{n_l})	\left[\varepsilon^h_b(E^h_{l},E^h_{k})+\delta_b(E_{l},E^h_{l})\right]^2.
            \end{split}
		\end{align}
	\end{lemma}

    \ENDO{\begin{remark}[On the Properties of the Error Bound]
The error bound in Lemma \ref{lem:eigenvec-estimation-algorithm-1} converges to zero as the approximate eigenvalues, $\lambda_{i,h}$, converge to the true ones, $\lambda_i$. The quality of the bound thus directly reflects the quality of the input eigenpairs, regardless of the approximation method.
While theoretically valid for any eigenvalue cluster, the bound is practically most effective for tight, well-separated clusters, as a wide cluster or a small spectral gap will yield a larger, less informative bound.
\end{remark}}
\medskip

        For the distances of subspaces, we have the following relation; see proof in Appendix.
	\begin{lemma}\label{lem:bbar-b-relation}
		For the distances $\delta_b$ and $\bar{\delta}_b$, it holds that
		\begin{gather}\label{eq:bbar-b-relation}
            \overline{\delta}^2_b(E_k,E^h_k)= 2-2\sqrt{1-\delta^2_b(E_k,E^h_k)}.
		\end{gather}
	\end{lemma}

        Additionally, we introduce the following lemma, which will be used to validate the linear independence of eigenfunctions:

        \begin{lemma}{[Lemma 2 of \cite{bronnimann1998interval}]}\label{lem:signal-of-det-estimation}
		Let $A$ be a $n\times n$ matrix. If $\|A-I\|<1$ for some matrix norm, then we have $\det(A)>0$, where $I$ denotes the $n\times n$ identity matrix.
	\end{lemma}

	\medskip
	% In addition, we will utilize the following lemma to validate the linear independence of eigenfunctions:
	% \begin{lemma}{[Lemma 2 of \cite{bronnimann1998interval}]}\label{lem:signal-of-det-estimation}
	% 	Let $A$ be a $n\times n$ matrix. If $\|A-I\|<1$ for some matrix norm, then we have $\det(A)>0$, where $I$ denotes the $n\times n$ identity matrix.
	% \end{lemma}
	% \medskip
								
	\section{Main theories}\label{section:main-theories}
								
        In this section, we focus on the case of multiple eigenvalues, i.e., $\lambda_n=\lambda_{n+1}=\cdots = \lambda_N$, and derive the difference quotient formula for eigenvalues, and provide an error estimate for the numerical approximations to the difference quotient formula. These results are used in Section 5 to evaluate the difference quotients of eigenvalues on nearly equilateral triangles, which leads to a partial solution to Conjecture \ref{len:main-conjecture}.
								
	\subsection{difference quotient formula of eigenvalues over triangles}
								
	For parameter $p=(x,y)(\in\mathbb{R}^2)$, let $T^p$ be a triangular domain with vertices as
	$O(0,0)$, $A(1,0)$ and $B(x,y)$. Suppose that the eigenvalues $\lambda_{n}^p,\cdots,\lambda_{N}^p$ are multiple at $p$, i.e., $\lambda_{n}^p=\cdots=\lambda_{N}^p(=:\lambda)$.								
	Let $E$ be the eigenspace corresponding to $\lambda$.
								
	For $t>0$ and $e\in\mathbb{R}^2$ with $\|e\|_2=1$, let $p_t:=p+te(=:(\tilde x,\tilde y))$ be a perturbation of $p$. To simplify the notation, let $T^t := T^{p_t}$, $T^0 := T^p$ and $\lambda_i^t:=\lambda_i^{p_t}$. Let $u_{n}^{t},\cdots,u_{N}^{t}$ be linearly independent eigenfunctions corresponding to the eigenvalues $\lambda_{n}^{t},\cdots,\lambda_{N}^{t}$, respectively.
    That is,
    \begin{equation}\label{eq:relation-lam-ut}
        \lambda_{i}^{t}= \frac{\|\nabla u_{i}^{t}\|^2_{T^t}}{\|u_{i}^{t}\|^2_{T^t}}~~(i=n,\cdots,N).
    \end{equation}
    Let us consider the linear transformation $S_t:\mathbb{R}^2\to\mathbb{R}^2$ that maps $T^0$ to $T^t$:
 	\begin{equation}
		S_t:=
		\begin{pmatrix}
			1 & (\tilde x-x)/y \\
			0 & \tilde y/y     
		\end{pmatrix}.
	\end{equation}
Let $\tilde u_{i}^t:=u^{t}_{i}\circ S_t (\in H^1_0(T^0))$ $(i=n,\cdots,N)$ and  define \begin{equation}\label{eq:def-widetildeE}
        \widetilde{E}_t:=\mbox{span}\{\tilde u_n^t,\cdots,\tilde u_N^t\} \quad  (\subset H_0^1(T^0)).
    \end{equation}
        Let
        \begin{equation}\label{eq:relation-lam-tilde-ut}
        \tilde\lambda_N^t:=\max_{\tilde u\in \widetilde E_t}\frac{\|\nabla\tilde u\|_{T^0}^2}{\|\tilde u\|_{T^0}^2}.
        \end{equation}

    Let us introduce several quantities to be used to characterize the difference quotient of eigenvalues 
    \begin{equation}
     D_t\lambda_i:=(\lambda_i^t-\lambda)/t~~(i=n,\cdots,N).
    \end{equation}
	Define $2\times 2$ matrices $P_t^e$ and $P^e$ by
	\begin{equation}\label{eq:def-PetPe}
		P_t^e:=\left(S_t^{-1}S_t^{-\intercal}-I\right)/t,~~P^e:=\lim_{t\to 0+}P^e_t.
	\end{equation}
	For a $2\times 2$ symmetric matrix $P$, introduce the symmetric bilinear form $F_{P}:H^1_0(T^0)\times H^1_0(T^0)\to\mathbb{R}$ defined by
	\begin{equation}\label{eq:def-FtPe}
		F_{P}(u,v):=(P\nabla u,\nabla v)_{T^0}~~\text{for } u,v\in H^1_0(T^0).
	\end{equation}
									
	\begin{remark}
		Matrix $P^e_t$ has the following expression with $e = (a, b)$:
		\begin{equation*}
			P_t^e = 
			\begin{pmatrix}
                \displaystyle\frac{ta^2}{\tilde{y}^2} & \displaystyle-\frac{ay}{\tilde{y}^2} \\
                \displaystyle-\frac{ay}{\tilde{y}^2} & \displaystyle-\frac{b(y + \tilde{y})}{\tilde{y}^2}
                \end{pmatrix}
		\end{equation*}
		Note that $(\tilde x,\tilde y)\to (x,y)$ $(t\to 0)$. Consequently, the matrix $P^e$ is given by
            \begin{equation}
            P^e=
            \begin{pmatrix}
            0 & \displaystyle-\frac{a}{y} \\
            \displaystyle-\frac{a}{y} & \displaystyle-\frac{2b}{y}
            \end{pmatrix}.
            \end{equation}		
	\end{remark}

        Next, we obtain the difference quotient formula of the eigenvalues.								 
	\begin{lemma}\label{lem:derivative-esimation-close-eigenvalues-prep}
		 $(D_t\lambda_i,\tilde u_{i}^t)$ ($i=n,\cdots,N$) are the eigenpairs of the following eigenvalue problem:

            Find $\tilde u\in \widetilde E_t$ and $\mu\in\mathbb{R}$ such that
		\begin{equation}\label{eq:derivative-esimation-close-eigenvalues-prep}
% before revision:
% (F_{P^{e}_{t}}(\tilde u,u)=) 
   (P^{e}_{t} \nabla \tilde u,\nabla u)_{T^0} 
			=
			\mu
			(\tilde u_,u)_{T^0}~~\mbox{$\forall u\in E$}.
		\end{equation}
	\end{lemma}
								
	\begin{proof}
		For each $i=n,\cdots,N$, the following variational equation holds:
		\begin{equation*}
			(\nabla u^t_{i},\nabla \tilde v)_{T^t}=\lambda^t_{i}(u^t_{i},\tilde v)_{T^t} \hspace{1em} \forall \tilde v\in H^1_0(T^t).
		\end{equation*}
		% that is, 
		% \begin{equation*}
		% (\widetilde\nabla u^{\tilde p},\widetilde\nabla(v\circ\Phi_{ p,\tilde p}^{-1}))_{T^t}=\lambda^{\tilde p}_k(u^{\tilde p},(v\circ\Phi_{ p,\tilde p}^{-1}))_{T^t} \hspace{1em} \forall v\in V(T^{ p}).
		% \end{equation*}
		Noting that $u^{t}_{i}=\tilde u_{i}^t\circ S_t^{-1}$ and
		$\nabla u^t_{i} = (S^{-\intercal}_{t})\nabla \tilde u_{i}^t$, we have
		\begin{equation}\label{eq:sub-above}
			\left(
			(S^{-\intercal}_{ t})\nabla\tilde u_{i}^t
			,
			(S^{-\intercal}_{t})\nabla v
			\right)_{T^0}\ENDO{\cdot |\det S_t|} 
			=
			\lambda^t_{i}
			(\tilde u_{i}^t,v)_{T^0}\ENDO{\cdot |\det S_t|}.
			\hspace{0.5em}\forall v\in H^1_0(T^0)~.
		\end{equation}
		For any $u\in E$, from \eqref{eq:sub-above}, it follows that
		\begin{equation*}
			\left(
			S^{-1}_{t}S^{-\intercal}_{t}\nabla\tilde u_{i}^t
			,
			\nabla u
			\right)_{T^0}\ENDO{\cdot |\det S_t|}
			=
			\lambda^t_{i}
			(\tilde u_{i}^t,u)_{T^0}\ENDO{\cdot |\det S_t|}.
		\end{equation*}
        \ENDO{Dividing by $|\det S_t|$, we obtain}
        \begin{equation}\label{eq:perturbated-variational-formula}
			\left(
			S^{-1}_{t}S^{-\intercal}_{t}\nabla\tilde u_{i}^t
			,
			\nabla u
			\right)_{T^0}
			=
			\lambda^t_{i}
			(\tilde u_{i}^t,u)_{T^0}.
		\end{equation}
        
		% Therefore, by taking the sum of the above equation for $k=n,\cdots,N$, we have
		% \begin{equation}
		%     \label{perturbated-formula}
		%     \sum_{k=n}^{N}s_{kj}\left(
		%     S^{-1}_{ p,\tilde p}S^{-\intercal}_{ p,\tilde p}\nabla\tilde u
		%     ,
		% \nabla u_{j}
		%     \right)_{T^0}
		%     =
		%     \sum_{k=n}^{N}
		%     s_{kj}\lambda^{\tilde p}_{k}
		%     (\tilde u,u_{j})_{T^0}.
		% \end{equation}
		Since $(\lambda, u)$ is an eigenpair over $T^0$, it holds that
		\begin{equation*}
			(\nabla u,\nabla v)_{T^0}=\lambda(u,v)_{T^0} \hspace{1em} \forall v\in H^1_0(T^0).
		\end{equation*}
		Taking $v=\tilde u^t_i(\in H^1_0(T^0))$ in the above variational equation,  we have
		\begin{equation}\label{eq:variational-formula}
			(\nabla u,\nabla\tilde  u_{i}^t)_{T^0}=\lambda(u,\tilde u_{i}^t)_{T^0}.
		\end{equation}
		Recall the definition of the symmetric bilinear form $F_{P^e_t}$ in \eqref{eq:def-FtPe}. From \eqref{eq:perturbated-variational-formula} and \eqref{eq:variational-formula}, we obtain
		\begin{equation}
			\left(
			\frac{1}{t}(S^{-1}_{t}S^{-\intercal}_{t}-I)\nabla\tilde u_{i}^t
			,
			\nabla u
			\right)_{T^0}
			=
			\frac{\lambda^t_{i}-\lambda}{t}
			(\tilde u_{i}^t,u)_{T^0},
		\end{equation}
		i.e.,
		\begin{equation}
			(P_t^e\nabla\tilde u_{i}^t,\nabla u)
			=
			D_t\lambda_i\cdot
			(\tilde u_{i}^t,u)_{T^0}~~\forall u\in E.
		\end{equation}
        Therefore, $(D_t\lambda_i,\tilde u_{i}^t)$ ($i=n,\cdots,N$) are the eigenpairs of the eigenvalue problem \eqref{eq:derivative-esimation-close-eigenvalues-prep}.
	\end{proof}

        \begin{remark}
            For the general case $\lambda_n^p\leq\cdots\leq\lambda_N^p$, rather than the eigenvalue problem \eqref{eq:derivative-esimation-close-eigenvalues-prep}, we can only obtain
            \begin{equation}
                F_{P^{e}_{t}}(\tilde u_{i}^t,u_{i})
			=
			D_t\lambda_i\cdot
			(\tilde u_{i}^t,u_{i})_{T^0}~(i=n,\cdots,N),
            \end{equation}
            where $u_{i}$'s are eigenfunctions corresponding to $\lambda^p_i~(i=n,\cdots,N)$.
        \end{remark}

                By taking concrete bases of the subspaces $E$ and $\widetilde E_t$, the variational equation for the difference quotients \eqref{eq:derivative-esimation-close-eigenvalues-prep} can be formulated as a matrix eigenvalue problem.
					
	\begin{theorem}\label{lem:Fte-basis-eigen}
		Suppose $\dim\widetilde E_t =\dim E$. 
		For a basis of $\widetilde E_t$, denoted by $\{\tilde\phi_i\}_{i=n}^{i=N}$, and a basis of $E$, denoted by $\{\phi_i\}_{i=n}^{i=N}$, let 
		\begin{equation}\label{eq:def-Mt-Nt}
			M_t:=\left(F_{P^{e}_{t}}\left(\tilde \phi_i,\phi_j\right)\right),~~
			N_t:=\left(\left( \tilde\phi_i,\phi_j\right)_{T^0}\right)~~(i,j=n,\cdots,N).
		\end{equation}
		Then, the difference quotient $D_t\lambda_i$ is the $(i-n+1)$-th eigenvalue of the following matrix eigenvalue problem:
		\begin{equation}\label{eq:Mt-N-eig-prob}
			M_t\sigma=\mu N_t\sigma.
		\end{equation}
							
		Moreover, for eigenfunctions $u_i^t~~(i=n,\cdots,N)$ over the perturbated domain $T^t$, we have
		\begin{equation*}
                \widetilde{u}_i^t=
			u_i^t\circ S_t=s_{ni}\tilde \phi_n+\cdots+s_{Ni}\tilde \phi_N~(\in\widetilde E_t)~~(i=n,\cdots,N),
		\end{equation*}
		where $\sigma_i:=(s_{ni},\cdots,s_{Ni})^\intercal$ becomes an eigenvector corresponding to the $(i-n+1)$-th eigenvalue of \eqref{eq:Mt-N-eig-prob}.
	\end{theorem}
	\begin{proof}
		From Lemma \ref{lem:derivative-esimation-close-eigenvalues-prep}, it follows that
		\begin{equation}\label{eq:F-eigen-phi}
			F_{P^{e}_{t}}(\tilde u_{i}^t,\phi_j)
			=
			D_t\lambda_i
			(\tilde u_{i}^t,\phi_j)_{T^0}~~(i,j=n,\cdots,N).
		\end{equation}              
		From the linear independence of $\{\tilde\phi_i\}_{i=n}^{i=N}$, each $\tilde u_i^t~~(i=n,\cdots,N)$ is uniquely expressed as
		\begin{equation}\label{tilde-u-t-tilde-phi}
			\tilde u_i^t=s_{ni}\tilde \phi_n+\cdots+s_{Ni}\tilde \phi_N,~~\mbox{where}~~(s_{ni},\cdots,s_{Ni})\in\mathbb{R}^{N-n+1}.
		\end{equation}
		Let $\sigma_i:=(s_{ni},\cdots,s_{Ni})^\intercal$.
		From the assumption $\dim\widetilde E_t =\dim E$, the system  $\{\sigma_i\}_{i=n}^{i=N}$ is linearly independent.
		Substituting \eqref{tilde-u-t-tilde-phi} into \eqref{eq:F-eigen-phi}, we obtain
		\begin{equation*}
			M_t\sigma_i=(D_t\lambda_i)N_t\sigma_i~~(i=n,\cdots,N).
		\end{equation*}
		It is clear that $D_e^{t}\lambda_n\leq D_e^{t}\lambda_{n+1}\leq,\cdots,\leq D_e^{t}\lambda_N$.
		From linear independence of  $\sigma_n,\cdots,\sigma_N$, the difference quotient $D_t\lambda_i~~(i=n,\cdots,N)$ is the $(i-n+1)$-th eigenvalue of $M_t\sigma=\mu N_t\sigma$, respectively.
	\end{proof}

        \medskip\medskip
					
	With further arguments on the differentiability of the eigenvalues, one can draw the same result of directional derivatives of simple or multiple eigenvalues as in Theorem \ref{lem:rousselet}.

\subsection{difference quotient formula for general polygons}

 In the case of general polygons, we can obtain analogous results about the difference quotients. For parameter $p=(x_1,x_2,\cdots,x_{k},y_1,y_2,\cdots,y_{k})\in\mathbb{R}^{2k}$, let $K^p$ be an $k$-gon with vertices $(x_i,y_i)~(i=1,\cdots,k)$.
				For $t>0$ and an $l^2$-normalized vector $e\in\mathbb{R}^{2k}$, let $p_t:=p+te$ be a perturbation of $p$. For simplicity of notations, denote by $K^t$  the polygonal domain $K^{p_t}$.
Let $K^{t}_h:=T_{1,t}\cup T_{2,t}\cup\cdots \cup T_{m,t}$ be a subdivision of $K^t$ by  disjoint triangles such that
    \begin{equation*}
        i\neq j\Rightarrow\mbox{int }T_{i,t}\cap\mbox{int }T_{j,t}=\emptyset.
    \end{equation*}
% For simplicity of notations, denote by $K^t$  the polygonal domain $K^{p_t}$, and by $T_{i,t}$ the triangular domain $T_i^{p_t}$ $(i=1,\cdots,n-2)$. 
Let $S_{j,t}$ be the linear transformation that maps the triangle $T_{j,0}$ to $T_{j,t}$ for each $j=1,\cdots,m$.
Let $\Phi_t:K^0\to K^t$ be the transformation that satisfies $\Phi_t|_{T_{j,0}}=S_{j,t}|_{T_{j,0}}$ for $j=1,\cdots,m$.
\ENDO{
We assume that the piecewise affine maps $\{S_{j,t}\}$ are defined such that the global transformation $\Phi_t$ is continuous across the boundaries of the subtriangles $\{T_{i,t}\}_{i=1}^m$. This ensures that if $u \in H_0^1(K^t)$, then its pullback $u \circ \Phi_t$ belongs to $H_0^1(K^0)$.}

                    \medskip
                    
                    Let $\lambda_i^t$ be the $i$-th Dirichlet eigenvalue of Laplacian over $K^{t}$.				Denote the difference quotient of the $i$-th eigenvalue $\lambda_i^p$ at $p$ by
\begin{equation}\label{def:nabla-e-t}
					D_t\lambda_i:=(\lambda^{t}_i-\lambda_i^0)/t.
				\end{equation}
                    Suppose that the eigenvalues $\lambda_{n}^p,\cdots,\lambda_{N}^p$ are multiple at $p$, i.e., $\lambda_{n}^p=\cdots=\lambda_{N}^p(=:\lambda)$.
				Let $E$ be the eigenspace corresponding to $\lambda$.

                    \medskip
							
				Let $u_{n}^{t},\cdots,u_{N}^{t}$ be linearly independent eigenfunctions corresponding to the eigenvalues $\lambda_{n}^{t},\cdots,\lambda_{N}^{t}$, respectively. Introduce $\tilde u_{i}^t=u^{t}_{i}\circ\Phi_t (\in H^1_0(K^0))$ $(i=n,\cdots,N)$. Note that 
					$\nabla (u^t_{i}|_{T_{j,t}}) = (S^{-\intercal}_{j,t})\nabla (\tilde u_{i}^t|_{T_{j,0}})$ for $j=1,\cdots, m$.
							 
				\begin{lemma}
					For the difference quotient of eigenvalues $D_t\lambda_i$ ($i=n,\cdots,N$), we have
					\begin{equation}\label{eq:derivative-esimation-close-eigenvalues-prep-poly}
						\tilde a_t(\tilde u_{i}^t,u)
						=
						D_t\lambda_i\cdot
						\tilde b_t(\tilde u_{i}^t,u)~~\mbox{$\forall u\in E$}.
					\end{equation}
     Here, $\tilde a_t,\tilde b_t$ are the symmetric bilinear forms  over $H^1_0(K^0)$ defined by
\begin{equation}\label{eq:def-tilde_at}
					\tilde a_t(u,v)=\sum_{j=1}^{m}\left[(P_{j,t}\nabla u,\nabla v)_{T_{j,t}}\cdot |\det S_{j,t}|+(\nabla u,\nabla v)_{T_{j,t}} d_{j,t}-\lambda( u,v)_{T_{j,t}}d_{j,t}\right],
				\end{equation}
                    \begin{equation}\label{eq:def-tilde_bt}
                        \tilde b_t(u,v)=\sum_{j=1}^{m}\left\{(u, v)_{T_{j,t}} d_{j,t}\right\},
                    \end{equation}
                    where $P_{j,t}:=( S_{j,t}^{-1}S_{j,t}^{-\intercal}-I)/t$ is a $2\times 2$ matrix, and $d_{j,t}:=(|\det S_{j,t}|-1)/t$.
				\end{lemma}
							
				\begin{proof}
     %Let us calculate the difference quotient $(\lambda^t_{i}-\lambda)/t$ ($i=n,\cdots,N$).
					% Introduce $\tilde u_{n(i)}^n=u^m_{n(i)}\circ\Phi_{ p,p_m}^{-1} \in H^1_0(T^{p_m})$.
					For each $i=n,\cdots,N$, the following variational equation holds:
					\begin{equation}\label{eq:perturbated-Kt-eig}
						(\nabla u^t_{i},\nabla \tilde v)_{K^t}=\lambda^t_{i}(u^t_{i},\tilde v)_{K^t} \hspace{1em} \forall \tilde v\in H^1_0(K^t).
					\end{equation}
					% that is, 
					% \begin{equation*}
					% (\widetilde\nabla u^{\tilde p},\widetilde\nabla(v\circ\Phi_{ p,\tilde p}^{-1}))_{K^t}=\lambda^{\tilde p}_k(u^{\tilde p},(v\circ\Phi_{ p,\tilde p}^{-1}))_{K^t} \hspace{1em} \forall v\in V(T^{ p}).
					% \end{equation*}
					 By changing the domain $K^t$ in \eqref{eq:perturbated-Kt-eig} to $K^0$, we have
					\begin{equation}\label{eq:sub-above-poly}
						\sum_{j=1}^{m}\left\{\left(
						(S^{-\intercal}_{j, t})\nabla\tilde u_{i}^t
						,
						(S^{-\intercal}_{j,t})\nabla v
						\right)_{T_{j,0}}\cdot|\det S_{j,t}|\right\}
						=
						\lambda^t_{i}
                            \sum_{j=1}^{m}\left\{
						(\tilde u_{i}^t,v)_{T_{j,0}}\cdot|\det S_{j,t}|\right\}
						\hspace{0.5em}\forall v\in H^1_0(K^0)~.
					\end{equation}
					For a given $u\in E$, substituting $v=u$ into \eqref{eq:sub-above-poly}, it follows that
					\begin{equation}\label{eq:perturbated-substitute-u-Kt}
						\sum_{j=1}^{m}\left\{\left(
						(S^{-1}_{j,t})(S^{-\intercal}_{j, t})\nabla\tilde u_{i}^t
						,
						\nabla u
						\right)_{T_{j,0}}\cdot|\det S_{j,t}|\right\}
						=
						\lambda^t_{i}
                            \sum_{j=1}^{m}\left\{
						(\tilde u_{i}^t,u)_{T_{j,0}}\cdot|\det S_{j,t}|\right\}
						.
					\end{equation}
					% Therefore, by taking the sum of the above equation for $k=n,\cdots,N$, we have
					% \begin{equation}
					%     \label{perturbated-formula}
					%     \sum_{k=n}^{N}s_{kj}\left(
					%     S^{-1}_{ p,\tilde p}S^{-\intercal}_{ p,\tilde p}\nabla\tilde u
					%     ,
					% \nabla u_{j}
					%     \right)_{K^0}
					%     =
					%     \sum_{k=n}^{N}
					%     s_{kj}\lambda^{\tilde p}_{k}
					%     (\tilde u,u_{j})_{K^0}.
					% \end{equation}
					Also, for the eigenpair $(\lambda, u)$ under the current setting, it holds that
					\begin{equation*}
						(\nabla u,\nabla v)_{K^0}=\lambda(u,v)_{K^0} \hspace{1em} \forall v\in H^1_0(K^0).
					\end{equation*}
					Taking $v=\tilde u^t_i$ in the above variational equation,  we have
					\begin{equation}\label{eq:variational-formula-poly}
						(\nabla u,\nabla\tilde  u_{i}^t)_{K^0}=\lambda(u,\tilde u_{i}^t)_{K^0}.
					\end{equation}
                    By substracting \eqref{eq:variational-formula-poly} from  \eqref{eq:perturbated-substitute-u-Kt} and then dividing by $t$, the left-hand side becomes
                    \begin{align}
                        \begin{split}\label{eq:substract-left}
                        \sum_{j=1}^{m}\frac{1}{t}&\left\{\left(\left(S_{j,t}^{-1}\right)\left(S_{j,t}^{-T}\right)\nabla\tilde u_i^t,\nabla u\right)_{T_{j,0}}\cdot|\det S_{j,t}|-\left(\nabla\tilde u_i^t,\nabla u\right)_{T_{j,0}}\right\}\\
                        &=
                        \sum_{j=1}^{m}\frac{1}{t}\left\{\left(\left(S_{j,t}^{-1}\right)\left(S_{j,t}^{-T}\right)\nabla\tilde u_i^t,\nabla u\right)_{T_{j,0}}\cdot|\det S_{j,t}|-\left(\nabla\tilde u_i^t,\nabla u\right)_{T_{j,0}}\cdot|\det S_{j,t}|\right\}\\
                        &~~+
                        \sum_{j=1}^{m}\frac{1}{t}\left\{\left(\nabla\tilde u_i^t,\nabla u\right)_{T_{j,0}}\cdot|\det S_{j,t}|-\left(\nabla\tilde u_i^t,\nabla u\right)_{T_{j,0}}\right\}\\
                        &=
                        \sum_{j=1}^{m}\left[\left(P_{j,t}\nabla\tilde u_i^t,\nabla u\right)_{T_{j,0}}\cdot|\det S_{j,t}|+\left(\nabla\tilde u_i^t,\nabla u\right)_{T_{j,0}}\cdot d_{j,t}\right].
                        \end{split}
                    \end{align}
                    The right-hand side becomes
                    \begin{align}
                        \begin{split}\label{eq:substract-right}
                        \sum_{j=1}^{m}&\frac{1}{t}\left\{\lambda^t_i(\tilde u_i^t,u)_{T_{j,0}}\cdot|\det S_{j,t}|-\lambda(\tilde u_i^t,u)_{T_{j,0}}\right\}\\
                        &=
                        \sum_{j=1}^{m}\frac{1}{t}\left\{\lambda^t_i(\tilde u_i^t,u)_{T_{j,0}}\cdot|\det S_{j,t}|-\lambda(\tilde u_i^t,u)_{T_{j,0}}\cdot|\det S_{j,t}|+\lambda(\tilde u_i^t,u)_{T_{j,0}}\cdot|\det S_{j,t}|-\lambda(\tilde u_i^t,u)_{T_{j,0}}\right\}\\
                        &=
                        \frac{\lambda^t_i–\lambda}{t}\sum_{j=1}^{m}\left\{(\tilde u_i^t,u)_{T_{j,0}}\cdot|\det S_{j,t}|\right\}+\sum_{j=1}^{m}\left\{\lambda(\tilde u_i^t,u)_{T_{j,0}}\cdot \frac{1}{t}\left(|\det S_{j,t}|-1\right)\right\}\\
                        \end{split}
                    \end{align}
					Recall that the symmetric bilinear form $F_t$ is defined by \eqref{eq:def-tilde_at} and \eqref{eq:def-tilde_bt}. From \eqref{eq:substract-left} and \eqref{eq:substract-right}, we obtain
					\begin{equation}
						\tilde a_t(\tilde u_{i}^t,u)
						=
						D_t\lambda_i\cdot
						\tilde b_t(\tilde u_{i}^t,u).
					\end{equation}
                     Therefore, $(D_t\lambda_i,\tilde u_{i}^t)$ ($i=n,\cdots,N$) are the eigenpairs of the eigenvalue problem \eqref{eq:derivative-esimation-close-eigenvalues-prep-poly}.
				\end{proof}

                \ENDO{\begin{remark}[On the dependency on $\lambda$]
The bilinear form $\tilde{a}_t$ in \eqref{eq:def-tilde_at} depends on the exact eigenvalue $\lambda$, which is typically unknown. For rigorous computation, one first computes a verified interval enclosure $[\underline{\lambda}, \overline{\lambda}]$ for $\lambda$. Substituting this interval into the definition of $\tilde{a}_t$ yields an interval matrix, reducing the problem to a generalized interval eigenvalue problem. Solving this with known algorithms provides a rigorous enclosure for the difference quotients $D_t\lambda_i$ that accounts for the initial uncertainty in $\lambda$.
\end{remark}}

        \medskip

	\subsection{Estimation for difference quotient formula over triangles}
        In this subsection, we provide an error estimate for the numerical approximations of the difference quotient formula in Theorem \ref{lem:Fte-basis-eigen}.
	We will continue to use the notation introduced in the previous subsection. 
					  
	Let $\hat u_n,\cdots,\hat u_N(\in H^1_0(T^0))$ be $L^2$-normalized approximations of the eigenfunctions corresponding to the multiple eigenvalues $\lambda_{n}^p=\cdots=\lambda_{N}^p(=\lambda)$. Let $\hat\lambda_{n}:=\min_{\hat u\in\widehat E}\|\nabla \hat u\|^2_{T_0}/\|\hat u\|_{T_0}^2$ $(i=n,\cdots,N)$ and \begin{equation}\label{eq:def-widehatE}
	    \widehat E:=\mbox{span}\{\hat u_n,\cdots,\hat u_N\}.
	\end{equation}
 Let $m:=N-n+1$. 
	Define $m\times m$ matrices $\widehat M_t$ and $\widehat N_t$ by
	\begin{equation}\label{def:approx-matrices}
		\widehat M_t=\left(F_{P_t^e}(\hat u_i,\hat u_j)\right)_{m\times m},~~\widehat N_t=\left((\hat u_i,\hat u_j)_{T^0}\right)_{m\times m}~~(i,j=n,\cdots,N).
	\end{equation}
	Note that $\widehat{M}_t,\widehat{N}_t$ will be constructed explicitly in the numerical computation.

To construct the matrices $M_t,N_t$ in the difference quotient formula \eqref{eq:Mt-N-eig-prob}, we choose specific systems $u_n^*,\cdots,u_N^*\in E$ and $\tilde u_n^*,\cdots,\tilde u_N^*\in \widetilde E_t$ as follow. 
        
        % To provide an error estimate for the difference quotients, we choose specific systems $u_n^*,\cdots,u_N^*\in E$ and $\tilde u_n^*,\cdots,\tilde u_N^*\in \widetilde E_t$ to construct the matrices $M_t,N_t$ in the difference quotient formula \eqref{eq:Mt-N-eig-prob}, which will be taken in the following way.    
        
\begin{description}
\item[Step 1:] Select $u_i^*\in E$ for each $\hat u_i\in \widehat E~(i=n,\cdots,N)$ such that
\begin{equation}\label{eq:b-estimation-ustar-uhat}
\|\hat u_i-u_i^*\|_{T^0}\leq \bar\delta_b(\widehat E,E),~~\|u_i^*\|_{T^0}=1.
\end{equation}
The definition of $\bar\delta_b$ ensures the existence of such elements $u_i^*$.

\item[Step 2:] Choose $\tilde u_i^*\in \widetilde E_t$ for each $u_i^*\in E~(i=n,\cdots,N)$ such that
\begin{equation}\label{eq:b-estimation-ustar-utilde}
\|u_i^*-\tilde u_i^*\|_{T^0}\leq \bar\delta_b(E,\widetilde E_t),~\|\tilde u_i^*\|_{T^0}=1.
\end{equation}
\end{description}
Note the linear independence of selected functions will be validated through Lemma \ref{lem:how-to-validate-the-linear-independence}. 

Let us summarize the relation among $\widehat E,E$ and $\widetilde E_t$ in Figure \ref{fig:relation-E-spaces}. Note that, in \eqref{eq:def-widehatE}, the space $\widehat E~(\subset H^1_0(T^0))$ is defined as a numerical approximation to the eigenspace $E~(\subset H^1_0(T^0))$, and in \eqref{eq:def-widetildeE}, the space $\widetilde E_t~(\subset H^1_0(T^0))$ is defined as the span of $\tilde u_i^t=u_i^t\circ S_t~(i=n,\cdots,N)$, where $u_i^t$ is an eigenfunction of $\lambda_i^t$ over the perturbated triangle $T^t$. 

\begin{figure}[H]
  \begin{minipage}[c]{0.45\linewidth}
    \centering
    \vskip 0.5cm
    \small 
    \begin{tikzpicture}[node distance=2cm]
 \node (E1) {$\widehat E$};
 \node (E2) [right=of E1] {$E$};
 \node (E3) [right=of E2] {$\widetilde E_t$};
 
 \node (u1) [below=0.5cm of E1] {$\hat u_i$};
 \node (u2) [below=0.5cm of E2] {$u_i^*$};
 \node (u3) [below=0.5cm of E3] {$\tilde u_i^*$};
 
 \node at ($(E1)!0.5!(u1)$) {\rotatebox{90}{$\in$}};
 \node at ($(E2)!0.5!(u2)$) {\rotatebox{90}{$\in$}};
 \node at ($(E3)!0.5!(u3)$) {\rotatebox{90}{$\in$}};
 
 \draw[->] (u1) -- (u2);
 \draw[->] (u2) -- (u3);
 
 \draw[-] (E1) to[out=30,in=150] node[midway,above] {$\overline{\delta}_b(\widehat E,E)$} (E2);
 \draw[-] (E2) to[out=30,in=150] node[midway,above] {$\overline{\delta}_b(E,\widetilde E_t)$} (E3);
\end{tikzpicture}
  \end{minipage}
  ~~
  \begin{minipage}[c]{0.49\linewidth}
    \centering 
    \footnotesize % または \footnotesize を使用
    \begin{align}
        \widehat E
        &=\mbox{span}\{\hat u_n,\cdots,\hat u_N\}\\
        E
        &=\mbox{span}\{u_n,\cdots,u_N\}
        =\mbox{span}\{u_n^*,\cdots,u_N^*\}\\
        \widetilde{E}_t
        &=\mbox{span}\{\tilde u_n^t,\cdots,\tilde u_N^t\}
        =\mbox{span}\{\tilde u_n^*,\cdots,\tilde u_N^*\}
    \end{align}
  \end{minipage}
  \caption{Relation among the spaces $\widehat E$, $E$, and $\widetilde E_t$.\label{fig:relation-E-spaces}}
\end{figure}            
The following matrices are used to evaluate \ENDO{the values of difference quotients}:
            \begin{equation}\label{eq:def-Mt-Nt-star}
			M_t^*:=\left(F_{P^{e}_{t}}\left(\tilde u_i^*,u_j^*\right)\right)_{m \times m},~~
			N_t^*:=\left(\left( \tilde u_i^*,u_j^*\right)_{T^0}\right)_{m \times m}~~(i,j=n,\cdots,N).
		\end{equation}
        Then, the eigenvalues of $M_t^*=\mu N_t^*$ are the same as the ones of \eqref{eq:Mt-N-eig-prob}, if the linear independence of $\{\tilde u_i^*\},\{u_i^*\}$ is validated; see the discussion in Lemma \ref{lem:how-to-validate-the-linear-independence}.
        We first provide error estimates for matrices $M_t^*$ and $N_t^*$, which allow us to evaluate the  difference quotients of eigenvalues.

        To bound the each element of $M^*_t$ and $N^*_t$, introduce the following quantities:
        \begin{align}
        \label{eq:def-EP}
        \begin{split}
        \hat\eta^2:=\lambda_N\overline{\delta}_b^2(\widehat E,E)+\hat{\lambda}_N-\lambda_n,~~
        \tilde\eta^2:=\lambda_N\overline{\delta}_b^2(E,\widetilde E_t)+\tilde{\lambda}_N^t-\lambda_n,
        \end{split}
        \end{align}
        where the quantity $\tilde{\lambda}_N^t$ is defined in \eqref{eq:relation-lam-tilde-ut}.
	\begin{lemma}\label{lem:M-star-N-star-estimation}
		For each element of matrices $M^*_t$ and $N^*_t$, we have
		\begin{gather}
			\begin{split}\label{ineq:F-b}
				\left|F_{P^e_t}(\tilde u_i^*,u_j^*)-F_{P^e_t}(\hat u_i,\hat u_j)\right|\leq \mbox{Err}_F(P^e_t;E,\widehat E)\\
				\text{ and }~\left|(\tilde u_i^*,u_j^*)_{T^0}-(\hat u_i,\hat u_j)_{T^0}\right|\leq \mbox{Err}_b(P^e_t;E,\widehat E),
			\end{split}
		\end{gather}
		where
		\begin{align}
			\mbox{Err}_F(P^e_t;E,\widehat E)&:=\|P^e_t\|_2\cdot(2\hat\eta\sqrt{\hat\lambda_{N}}+\tilde\eta\sqrt{\lambda_{N}}),\notag\\
            \mbox{Err}_b(P^e_t;E,\widehat E)&:=2\overline{\delta}_b(\widehat E,E)+\overline{\delta}_b(E,\widetilde E_t).
		\end{align}
	\end{lemma}
					
\begin{proof}
For each $i,j=n,\cdots,N$, from the estimation \eqref{eq:b-estimation-ustar-uhat}, we obtain
\begin{align}\label{eq:elementwise-N-I}
\begin{split}
&\biggl|
\left(\tilde u^*_{i},u_{j}^*\right)
-
\left(\hat u_{i},\hat u_{j}\right)
\biggr|
=    
\biggl|
\left(\tilde u^*_{i}, u_{j}^*\right)_{T^0}
-
\left(\hat u_{i},\hat u_{j}\right)_{T^0}
\biggr|\\
& \leq 
\biggl|
\left(\tilde u^*_{i}, u_{j}^*\right)_{T^0}
-
\left( u_{i}^*, u_{j}^*\right)_{T^0}
\biggr|+
\biggl|
\left( u_{i}^*, u_{j}^*\right)_{T^0}
-
\left(\hat u_{i}, u_{j}^*\right)_{T^0}
\biggr|+
\biggl|
\left(\hat u_{i}, u_{j}^*\right)_{T^0}
-
\left(\hat u_{i},\hat u_{j}\right)_{T^0}
\biggr|
\\
& \leq 
\left\|\tilde u_{i}^*- u_{i}^*\right\|_{T^0}\left\| u_{j}^*\right\|_{T^0}+
\left\| u_{i}^*-\hat u_{i}\right\|_{T^0}\left\| u_{j}^*\right\|_{T^0}+
\left\|\hat u_{i}\right\|_{T^0}\left\| u_{j}^*-\hat u_{j}\right\|_{T^0}\\
& \leq 
2\overline{\delta}_b(\widehat E,E)+\overline{\delta}_b(E,\widetilde E_t)(=\mbox{Err}_b(P^e_t;E,\widehat E)).
\end{split}
\end{align}
From the fundamental formula \eqref{eq:fundamental-formula}, for $i=n,\cdots,N$, we have
\begin{align}\label{eq:est-star-hat}
\left\|\nabla u_i^*-\nabla \hat{u}_i^*\right\|_{T^0}^2
&\leq\lambda_N\left\|u_i^*-\hat{u}_i^*\right\|_{T^0}^2-\left(\lambda_n-\hat{\lambda}_N\right)\left\|\hat{u}^*_i\right\|_{T^0}^2\\
&\leq
\lambda_N\overline{\delta}_b^2(\widehat E,E)+\hat{\lambda}_N-\lambda_n(=\hat\eta^2),
\end{align}
and
\begin{align}\label{eq:est-tilde-star}
\left\|\nabla\tilde u_i^*-\nabla u_i^*\right\|_{T^0}^2
&\leq\lambda_N\left\|\tilde u_i^*-u_i^*\right\|_{T^0}^2-\left(\lambda_n-\tilde\lambda_N^t\right)\left\|\hat{u}^*_i\right\|_{T^0}^2\\
&\leq
\lambda_N\overline{\delta}_b^2( E,\widetilde E_t)+\tilde\lambda_N^t-\lambda_n(=\tilde\eta^2).
\end{align}

Furthermore, from estimations \eqref{eq:est-star-hat} and \eqref{eq:est-tilde-star}, it holds that
\begin{align}\label{eq:elementwise-M-M-star}
\begin{split}
\biggl|
F_{P^e_t}\left(\tilde u^*_{i},u^*_{j}\right)
-
F_{P^e_t}\left(\hat u_{i},\hat u_{j}\right)
\biggr|
& =      
\biggl|
\left(P^e_t\nabla\tilde u^*_{i},\nabla u^*_{j}\right)_{T^0}
-
\left(P^e_t\nabla\hat u_{i},\nabla\hat u_{j}\right)_{T^0}
\biggr|\\
& \leq   
\biggl|
\left(P^e_t\nabla\tilde u^*_{i},\nabla u^*_{j}\right)_{T^0}
-
\left(P^e_t\nabla u^*_{i},\nabla u^*_{j}\right)_{T^0}
\biggr|\\
& \quad+ 
\biggl|
\left(P^e_t\nabla u^*_{i},\nabla u^*_{j}\right)_{T^0}
-
\left(P^e_t\nabla\hat u_{i},\nabla u^*_{j}\right)_{T^0}
\biggr|
\\
& \quad+ 
\biggl|
\left(P^e_t\nabla\hat u_{i},\nabla u^*_{j}\right)_{T^0}
-
\left(P^e_t\nabla\hat u_{i},\nabla\hat u_{j}\right)_{T^0}
\biggr|
\\
& \leq   
\|P^e_t\|_2\left\|\nabla\tilde u_{i}^*-\nabla u^*_{i}\right\|_{T^0}\left\|\nabla u^*_{j}\right\|_{T^0}\\
& \quad+ 
\|P^e_t\|_2\left\|\nabla u^*_{i}-\nabla\hat u_{i}\right\|_{T^0}\left\|\nabla u^*_{j}\right\|_{T^0}\\
& \quad+ 
\|P^e_t\|_2\left\|\nabla\hat u_{i}\right\|_{T^0}\left\|\nabla u^*_{j}-\nabla\hat u_{j}\right\|_{T^0}\\
& \leq   
\|P^e_t\|_2\cdot(2\hat\eta\sqrt{\hat\lambda_{N}}+\tilde\eta\sqrt{\lambda_{N}})~(=\mbox{Err}_F(P^e_t;E,\widehat E)).
\end{split}
\end{align}

\end{proof}
        To validate the linear independence of the systems $u_n^*,\cdots,u_N^*\in E$ and $\tilde u_n^*,\cdots,\tilde u_N^*\in \widetilde E_t$, we verify the non-singularity of the matrices $N^*_t$ and $\widetilde N^*_t$ defined by
        \begin{equation*}
        N^*_t:=\left(\left(u^*_{i},u^*_{j}\right)_{T^0}\right),~~
			\widetilde N^*_t:=\left(\left(\tilde u^*_{i},\tilde u^*_{j}\right)_{T^0}\right)~~(i,j=n,\cdots,N).
		\end{equation*}

    \ENDO{
	The following lemma is used to validate the linear independence of these systems.
\begin{lemma}\label{lem:how-to-validate-the-linear-independence}
    For the systems $u_n^*, \ldots, u_N^* \in E$ and $\tilde{u}_n^*, \ldots, \tilde{u}_N^* \in \widetilde{E}_t$, the following properties hold:
    \begin{description}
        \item[(i)] If $\mathrm{Est}(N_t^*) < 1$, then $u_n^*, \ldots, u_N^* \in E$ are linearly.
        \item[(ii)] If $\mathrm{Est}(\widetilde{N}_t^*) < 1$, then $\tilde{u}_n^*, \ldots, \tilde{u}_N^* \in \widetilde{E}_t$ are linearly independent.
    \end{description}
\end{lemma}
 \begin{proof}
            To prove (i), suppose $\mbox{Est}(\widehat N_t,N_t^*)<1$. 
            For each element of $\widehat N_t$ and $N^*_t$, from \eqref{eq:b-estimation-ustar-uhat}, we have
		\begin{align*}
			\left|
			\left(\hat u_{i},\hat u_{j}\right)_{T^0}
			-
			\left(u^*_{i}, u^*_{j}\right)_{T^0}
			\right|
			  & \leq   
			\left|(\hat u_{i},\hat u_{j})_{T^0}-(u_{i}^*,\hat u_{j})_{T^0}\right|
			+
			\left|(u_{i}^*,\hat u_{j})_{T^0}-(u_{i}^*,u_{j}^*)_{T^0}\right|
			\\
			  & \leq   
			\left\|\hat u_{i}-u_{i}^*\right\|_{T^0}\left\|\hat u_{j}\right\|_{T^0}
			+\left\|u_{i}^*\right\|_{T^0}\left\|\hat u_{j}-u_{j}^*\right\|_{T^0}\\
			  & \leq   
			2\overline{\delta}_b(\widehat E,E).
		\end{align*}
            Then, we have
            \begin{align}
                \left\|I- N_t^*\right\|_\infty  
                &\leq
                \left\|I-\widehat N_t\right\|_\infty+\left\|\widehat N_t-N_t^*\right\|_\infty\\
                &\leq
                \left\|I-\widehat N_t\right\|_\infty+2m\overline{\delta}_b(\widehat E,E)\\
                &=
                \mbox{Est}(\widehat N_t,N_t^*)<1
            \end{align}
            Therefore, by Lemmma \ref{lem:signal-of-det-estimation}, $u_n^*,\cdots,u_N^*\in E$ is linearly independent.
            
		To prove (ii), suppose $\mbox{Est}(\widehat N_t,\widetilde N_t^*)<1$.  For each element of $\widehat N_t$ and $\widetilde{N}^*_t$, by \eqref{eq:b-estimation-ustar-utilde}, we have
		\begin{align*}
			\left|
			\left(\hat u_{i},\hat u_{j}\right)
			-
			\left(\tilde u^*_{i},\tilde u^*_{j}\right)
			\right|
			  & =      
			\left|
			\left(\hat u_{i},\hat u_{j}\right)_{T^0}
			-
			\left(\tilde u^*_{i}, \tilde u^*_{j}\right)_{T^0}
			\right|\\
			  & \leq   
			\left|(\hat u_{i},\hat u_{j})-(u^*_{i}, u^*_{j})\right|
			+
			\left|(u_{i}^*, u_{j}^*)-(\tilde u_{i}^*,\tilde u_{j}^*)\right|
			\\
			  & \leq   
			\left\|\hat u_{i}-u_{i}^*\right\|_{T^0}\left\|\hat u_{j}\right\|_{T^0}
			+\left\|u_{i}^*\right\|_{T^0}\left\|\hat u_{j}-u_{j}^*\right\|_{T^0}\\
			  & \quad+ 
			\left\|u_{i}^*-\tilde u_{i}^*\right\|_{T^0}\left\|u_{j}^*\right\|_{T^0}
			+\left\|\tilde u_{i}^*\right\|_{T^0}\left\|u_{j}^*-\tilde u_{j}^*\right\|_{T^0}\\
			  & \leq   
			2(\overline{\delta}_b(\widehat E,E)+\overline{\delta}_b(E,\widetilde E_t)).
		\end{align*}
            Then, we have
            \begin{align}
                \left\|I- N_t^*\right\|_\infty  
                &\leq
                \left\|I-\widehat N_t\right\|_\infty+\left\|\widehat N_t-\widetilde N_t^*\right\|_\infty\\
                &\leq
                \left\|I-\widehat N_t\right\|_\infty+2m(\overline{\delta}_b(\widehat E,E)+\overline{\delta}_b(E,\widetilde E_t))\\
                &=
                \mbox{Est}(\widehat N_t,\widetilde N_t^*)<1
            \end{align}
            Therefore, from Lemma \ref{lem:signal-of-det-estimation}, $\tilde u_n^*,\cdots,\tilde u_N^*\in E$ is linearly independent.
	\end{proof}
    }
Table \ref{table:summary} summarizes the notations and their roles in the error estimation framework.

\begin{table}[h]
 \caption{Summary of notations for error estimation.\label{table:summary}}
 \centering
 \begin{tabular}{|c|p{5.5cm}|p{5.5cm}|}
  \hline
  Symbol & Description & Role and Relation \\
  \hline
  $\widehat{M}_t, \widehat{N}_t$ & Matrices constructed from numerical approximations of eigenfunctions, $\{\hat{u}_i\}$. & These are the matrices we can explicitly calculate. \\
  \hline
  $M_t^*, N_t^*$ &``True" matrices for the difference quotient eigenvalue problem, built from the ideal basis functions $\{u_i^*\}$ and $\{\tilde{u}_i^*\}$. & The objects we aim to find rigorous bounds for. \\
  \hline
  $\overline{\delta}_b(\widehat{E}, E)$ & Directed distance from the approximate eigenfunction space $\widehat{E}$ to the exact space $E$. & Measures the numerical approximation error of the eigenspace. \\
  \hline
  $\overline{\delta}_b(E, \widetilde{E}_t)$ & Directed distance from the exact space $E$ to the perturbed exact space $\widetilde{E}_t$. & Measures the domain perturbation error caused by changing the domain from $T^0$ to $T^t$. \\
  \hline
  $\mathrm{Err}_F, \mathrm{Err}_b$ & The final, computable error bounds for the matrix elements, defined in Lemma 3.4. & Provide the rigorous enclosure: \newline $|(M_t^*)_{ij} - (\widehat{M}_t)_{ij}| \le \mathrm{Err}_F(P^e_t;E,\widehat E)$, \newline $|(N_t^*)_{ij} - (\widehat{N}_t)_{ij}| \le \mathrm{Err}_b(P^e_t;E,\widehat E)$. \\
  \hline
 \end{tabular}
\end{table}

\vspace{15em}

\ENDO{
        \subsection{Estimation for the errors $\mbox{Err}_F$ and $\mbox{Err}_b$}
        In this section, we provide the estimates for the quantities that are required to be obtained to evaluate each element of the matrices $M^*_t$ and $N^*_t$ in Lemma \ref{lem:M-star-N-star-estimation}.

        In Lemma \ref{lem:M-star-N-star-estimation}, we introduced the following quantities to bound each element of the matrices $M^*_t$ and $N^*_t$:
        \begin{align}
            \mbox{Err}_F(P^e_t;E,\widehat E)&:=\|P^e_t\|_2\cdot(2\hat\eta\sqrt{\hat\lambda_{N}}+\tilde\eta\sqrt{\lambda_{N}}),\notag\\
            \mbox{Err}_b(P^e_t;E,\widehat E)&:=2\overline{\delta}_b(\widehat E,E)+\overline{\delta}_b(E,\widetilde E_t)\label{eq:def-Err-Fb}.
        \end{align}
        The quantities $\hat\eta,\tilde\eta$ are defined in  \eqref{eq:def-EP} as follows:
        \begin{align}
        \begin{split}
        \hat\eta^2:=\lambda_N\overline{\delta}_b^2(\widehat E,E)+\hat{\lambda}_N-\lambda_n,~~
        \tilde\eta^2:=\lambda_N\overline{\delta}_b^2(E,\widetilde E_t)+\tilde{\lambda}_N^t-\lambda_n.
        \end{split}
        \end{align}
        Here, $\tilde{\lambda}_N^t$ is the quantity defined in \eqref{eq:relation-lam-tilde-ut}, i.e.,
        \begin{equation*}
        \tilde\lambda_N^t:=\max_{\tilde u\in \widetilde E_t}\frac{\|\nabla\tilde u\|_{T^0}^2}{\|\tilde u\|_{T^0}^2}.
        \end{equation*}

        The quantity $\tilde{\lambda}_N^t$ can be estimated by using the following lemma:
        \begin{lemma}\label{eq:bound-tilde-lambda}
        For the quantity $\tilde{\lambda}_N^t$, we have the following upper bound:
        \begin{equation}\label{eq:bound-relation-lam-tilde-ut}
        \tilde\lambda_N^t\leq \|S_tS_t^\intercal\|^2_2\cdot\lambda_N^t,
        \end{equation}
        where $S_t:\mathbb{R}^2\to\mathbb{R}^2$ is the linear transformation that maps $T^0$ to $T^t$.
        \end{lemma}
        \begin{proof}
        Let $E^t:=\mbox{span}\{u_n^t\cdots u_N^t\}$. Then, we have
        \begin{align*}
            \tilde\lambda_i^t
            &=\max_{\tilde u\in\widetilde   E_t}\frac{\|\nabla\tilde u\|_{T^0}^2}{\|\tilde u\|_{T^0}^2}
            =\max_{u\in E^t}\frac{\|S_tS_t^\intercal\nabla u\|_{T^t}^2}{\| u\|_{T^t}^2}\\
            &\leq
            \|S_tS_t^\intercal\|^2_2\max_{u\in E^t}\frac{\|\nabla u\|_{T^t}^2}{\| u\|_{T^t}^2}
            =
            \|S_tS_t^\intercal\|^2_2\cdot\lambda^t_N.
        \end{align*}
        \end{proof}

        Upper bounds for $\hat\eta$ and $\tilde\eta$ are estimated using Lemma \ref{lem:eigenvec-estimation-algorithm-1} and \ref{lem:bbar-b-relation}. Note that for $\hat\eta$, the concrete approximate eigenfunction $\hat u_i$'s will be used to evaluate $\hat\lambda_i$'s. For $\tilde\eta$, since the basis of $\widetilde E_t$ cannot be computed directly, we show the idea to compute an upper bound of $\tilde\eta$.
        
        \medskip
        \textbf{Estimation of $\overline{\delta}_b(E,\widetilde E_t)$: }
        \medskip
        
        From the relation \eqref{eq:bbar-b-relation}, to estimate the upper bound of $\overline{\delta}_b(E,\widetilde E_t)$, it suffices to compute the upper bound of $\delta_b(E,\widetilde E_t)$. 
        The estimation of  $\delta_b(E,\widetilde E_t)$ is obtained by applying Lemma \ref{lem:eigenvec-estimation-algorithm-1} in a recursive way. To simplify the argument,  let us assume that $\lambda_1, \cdots, \lambda_{n-1}$ forms the first cluster of eigenvalues. That is, 
        \begin{align}
            E_1 &:= \mbox{span}\{u_1,u_2,\cdots,u_{n-1}\},~~
            E_2 := E,\\
            E^{h}_1 &:= \mbox{span}\{\tilde u_1^t,\tilde u_2^t,\cdots,\tilde u_{n-1}^t\},~~
            E^{h}_2 := \widetilde{E}_t.\\
        \end{align}
        Then, from Lemma \ref{lem:eigenvec-estimation-algorithm-1}, we have 
        \begin{equation}\label{eq:delta-b-e2-e2h}
            \delta^2_b(E_2,E^h_{2})\leq\frac{\tilde\lambda_N^t-\lambda_n+\theta_b}{\rho-\lambda_n},
        \end{equation}
        where $\rho>0$ is a quantity such that $\lambda_{n}<\rho\leq\lambda_{N+1}$, and $\theta_b$ is defined by
        \begin{equation}
            \theta_b=(\rho-\lambda_n)\left[\varepsilon^h_b(E^h_{1},E^h_2)+\delta_b(E_{1},E^h_1)\right]^2.
        \end{equation}
        From the definition of $\widetilde{u}_i^t$, it is easy to see that the spaces $E_1^h$ and $E_2^h$ are $L^2$-orthogonal to each other. Hence, 
        \begin{equation}
            \varepsilon^h_b(E^h_{1},E^h_2):=\max_{\substack{v\in E^h_{1}\\\|v\|_{T^0}=1}}\max_{\substack{v_h\in E^h_{2}\\ \|v_h\|_{T^0}=1}}(v,v_h)_{T^0}=0.
        \end{equation}
        
    }        
							
\section{Computer-assisted proof for the simplicity of $\lambda_2$}

	\label{section:computation-results}

                \subsection{Preparation}

        In this section,  by using the rigorous estimation for difference quotients provided in the previous section, we propose a computer-assisted proof for the following theorem:
        \medskip
        \begin{theorem}\label{len:partial-result}
        The second Dirichlet eigenvalue is simple for every non-equilateral triangle with its minimum normalized height greater than or equal to $\tan(\pi/60)/2$.
        \end{theorem}
	\medskip
    
\begin{figure}[H]
		\centering
		\includegraphics[keepaspectratio, scale=0.27]{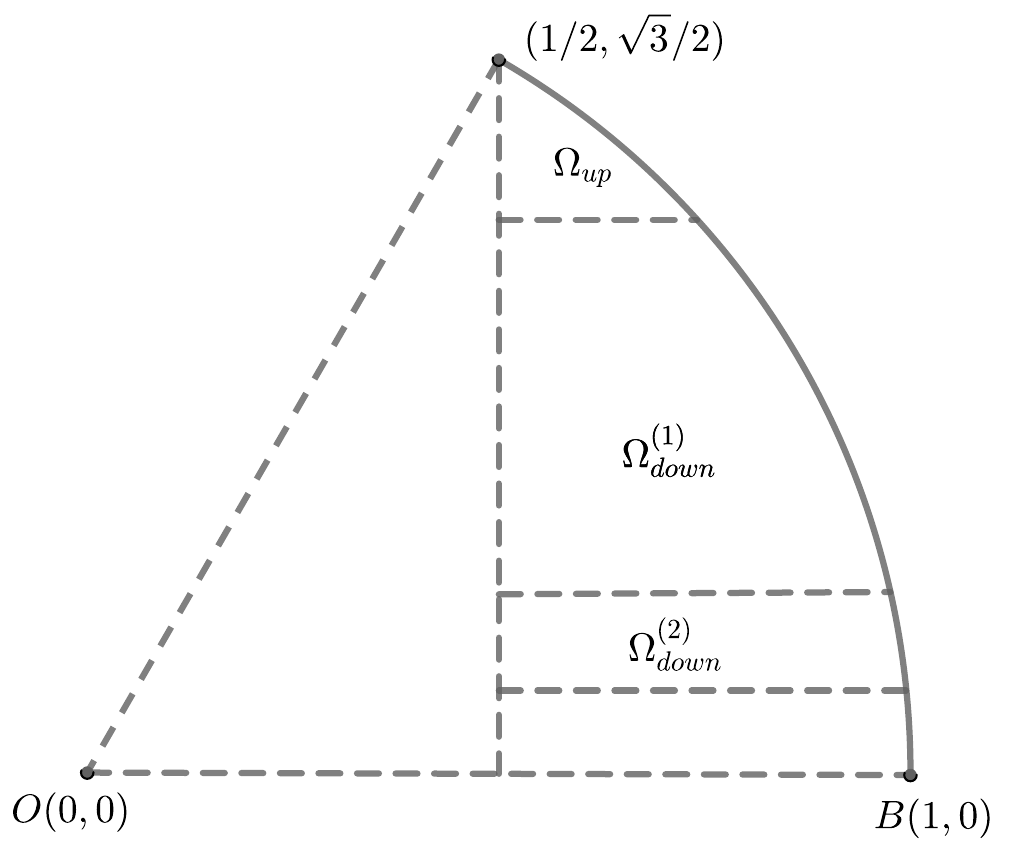} ~ 
  	\includegraphics[keepaspectratio, scale=0.45]{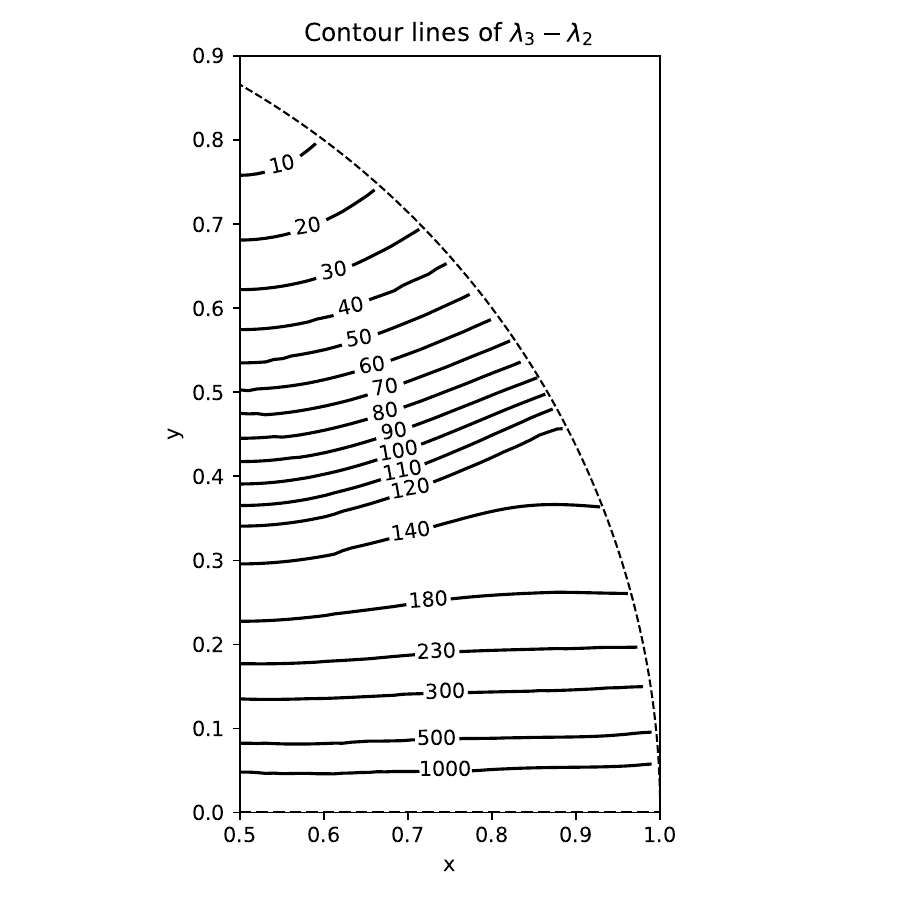}
		\caption{Subdivisions of $\Omega$ (left) and contour lines of $\lambda_3^p-\lambda_2^p$ (right) \label{fig:omega-up-down}}
	\end{figure}

	Since the simplicity of eigenvalues is isometry and scaling invariant, without loss of generality, the vertex $p$ of a triangle $T^p$ is further assumed to be located in the moduli space of triangles defined by
	\begin{equation*}
		\Omega:=\{(x,y)\in\mathbb{R}^2:x^2+y^2\leq 1, x\geq 1/2,y>0\}.
	\end{equation*}
        Let us subdivide $\Omega$ into four parts $\Omega_{up},\ENDO{\Omega_{down}^{(1)}},\ENDO{\Omega_{down}^{(2)}}$ and $\ENDO{\Omega_{rest}}$. 
	\begin{align*}
		\Omega_{up} :&= \{(x,y) \in \Omega ~:~ y \geq (\sqrt{3} - \varepsilon)/2\},\\
		\ENDO{\Omega_{down}^{(1)}} :&= \left\{(x,y) \in \Omega ~:~
		\tan\left(\frac{\pi}{14}\right)\leq y < (\sqrt{3} - \varepsilon)/2\right\},\\
        \ENDO{\Omega_{down}^{(2)}} :&= \left\{(x,y) \in \Omega ~:~
		\frac{1}{2}\tan\left(\frac{\pi}{60}\right)\leq y\leq\tan\left(\frac{\pi}{14}\right)\right\},\\
        \ENDO{\Omega_{rest}} :&= \left\{(x,y) \in \Omega ~:~ 
		y<\frac{1}{2}\tan\left(\frac{\pi}{60}\right)\right\}.
	\end{align*}
 	The value of $\varepsilon>0$ will be explicitly determined in the process of verified computation. For the geometric settings described above, see the left part of Figure \ref{fig:omega-up-down}.

In the rest of this section, we will investigate the simplicity of $\lambda_2^p$ over $\Omega_{up}$ and $\ENDO{\Omega_{down}^{(1)}}\cup \ENDO{\Omega_{down}^{(2)}}$ by using two different approaches. 
The simplicity of $\lambda_2$ can be observed in Figure \ref{fig:omega-up-down}, which displays the contour lines of the numerically computed value of  $\lambda_3^p -\lambda_2^p$.  

\medskip

        Let $p_0=(1/2,\sqrt{3}/2)$. $T^{p_0}$ is an equilateral triangle.
        To prove Theorem \ref{len:partial-result}, it suffices to show that  $\lambda_2^p<\lambda_3^p$ for every $p\in(\Omega_{up}\cup\ENDO{\Omega_{down}^{(1)}\cup \Omega_{down}^{(2)}})\backslash\{p_0\}$.

 \medskip
 
	Recall that the difference quotient of the $k$th eigenvalue is denoted by
	\begin{equation}\label{def:nabla-e-t-1}
		D\lambda_k(p,\tilde p):=(\lambda^{\tilde p}_k-\lambda_k^p)/\|\tilde p-p\|_2~~(p,\tilde p\in\Omega,~p\neq\tilde p).
	\end{equation}
	We provide a computer-assisted proof for Theorem \ref{len:partial-result} by following the outline below:
					
	\begin{description}
\item[Step 1:] (Case of $\Omega_{up}$) By concretely evaluating the difference quotient for the eigenvalues, prove that the difference quotients of the second and third eigenvalues satisfy
\begin{equation}
D\lambda_2(p_0,p) < D\lambda_3(p_0,p) \quad \text{for all } p \in \Omega_{up}\backslash\{p_0\},
\end{equation}
that is,
\begin{equation*}\label{eq:separation-diff-quot}
    \frac{\lambda_2^p-\lambda_2^{p_0}}{\|p-p_0\|_2}<\frac{\lambda_3^p-\lambda_3^{p_0}}{\|p-p_0\|_2} \quad \text{for all } p \in \Omega_{up}\backslash\{p_0\}.
\end{equation*}
Since $\lambda_2^{p_0} = \lambda_3^{p_0}$ for the equilateral triangle $T^{p_0}$ with $p_0=(1/2,\sqrt{3}/2)$, this implies $\lambda_2^{p} < \lambda_3^{p}$ for all $p \in \Omega_{up}\backslash{p_0}$.
\item[Step 2:] (Case of $\ENDO{\Omega_{down}^{(1)}}\cup \ENDO{\Omega_{down}^{(2)}}$) Separate $\lambda_2^p$ and $\lambda_3^p$ by directly estimating the upper and lower bounds of the eigenvalues:
\begin{equation}
\lambda_2^p\leq \overline{\lambda}_2^p <\underline{\lambda}_3^p\leq \lambda_3^p \quad \text{for all } p \in \ENDO{\Omega_{down}^{(1)}}\cup \ENDO{\Omega_{down}^{(2)}}.
\end{equation}
The eigenvalue bounds $ \overline{\lambda}_2^p$ and $ \underline{\lambda}_3^p$ will be calculated directly through verified computation.
\end{description}

\medskip

\paragraph{Approximation of eigenpairs by Finite Element Method (FEM)}
The FEM will be utilized to  approximate the eigenvalues and their corresponding eigenfunctions over triangles. 
For a triangular domain $T\subset\mathbb{R}^2$, denote by $\mathcal{T}^h$ a regular triangulation of $T$; that is, any two edges $e_i$ and $e_j$ of elements of $\mathcal{T}^h$ satisfy $e_i\cap e_j=e_i=e_j$ or $\mu(e_i\cap e_j)=0$, where $\mu(\cdot)$ is the $1$-dimensional measure.  Let $h$ be the maximal edge length of $\mathcal{T}^h$.
The space of polynomials of degree up to $d$ over an element $K$ of $\mathcal{T}_h$ is denoted by $P^d(K)$.
                    
Let us introduce the following finite element space $V_h^{\CG}$ over $\mathcal{T}^h$:
\begin{itemize}
\item The Lagrange FEM space $V_h^{\CG,d}(T)$ of degree $d(\ge 1)$:
\end{itemize}
\begin{equation}
\label{def:fem-space-cg}
V_h^{\CG,d}(T):=
\{\hat v:\hat v\mbox{ is a continuous piecewise polynomial of degree $d$ on } \mathcal{T}^h.
\}
  \end{equation}
Also, define $V_{h,0}^{\CG,d}(T) := H^1_0(T) \cap V_{h}^{\CG,d}(T)$. For simplicity, the notation will be abbreviated as $V_{h}^{\CG}$ or $V_{h,0}^{\CG}$.

To approximate the eigenvalues and the eigenfunctions of \eqref{eq:eigenvalue-problem}, we consider the following discretized eigenvalue problem:
                    
\begin{itemize}\label{eq:lambda-h}
\item Find $\hat u\in V^{\CG}_{h,0}(T)\backslash\{0\}$ and $\hat \lambda>0$  such that
      \begin{equation}\label{eq:cg-problem}
        (\nabla \hat u,\nabla \hat v)_{T}=\hat \lambda (\hat u,\hat v)_{T} \hspace{1em} \forall \hat v\in V^{\CG}_{h,0}(T).
      \end{equation}
\end{itemize}
Let  $N_{1}=\mbox{dim}(V_{h,0}^{\CG}(T))$.
The eigenvalues of the above problem are denoted by
\begin{equation*}
(0<)~\hat \lambda_{1}^{\CG}(T)\leq\hat \lambda_{2}^{\CG}(T) \leq \cdots \leq \hat\lambda_{N_1}^{\CG}(T)~.
\end{equation*}
Let $\hat u_1,\hat u_2,\cdots,\hat u_{N_1}\in V_{h,0}^\CG$ be eigenfunctions corresponding to $\hat \lambda_{1}^{\CG}(T), \cdots ,\hat\lambda_{N_1}^{\CG}(T)$, respectively.

\medskip

\begin{lemma}[\ENDO{Theorem 4.3 of \cite{liu2013verified}}]\label{lem:est-tau}
    For triangulation $\mathcal{T}^h$ consisting of right triangles, we have
    \begin{equation}\label{lem:l-estimation-cg}
\underline\lambda_k:=\frac{\hat\lambda_{k}^{\CG}(T)}{1+(0.493h)^2\hat\lambda_{k}^{\CG}(T)}\leq\lambda_k(T)\leq\hat\lambda_{k}^{\CG}(T)=:\overline\lambda_k        
        \hspace{1em} 
        \mbox{ for } k=1,2,...,N_1,
    \end{equation}
where $h$ is the mesh size of $\mathcal{T}^h$; see Figure \ref{fig:right-triangular-mesh} for the specific mesh shape. 
\end{lemma}

\paragraph{High-precision eigenvalue bounds by Lehmann-Goerisch's Method}

Lemma \ref{lem:est-tau} can provide rigorous lower eigenvalue bounds.  However, this result cannot take advantage of high-degree finite element methods, and one has to refine the mesh to improve the precision of the computed lower eigenvalue bounds.  

To obtain highly precise eigenvalue bounds in an efficient way, we follow the approach proposed in \cite{liu2024lehmann}, which utilizes the Lehmann-Goerisch method \cite{goerisch1990determination,behnke1994inclusions,GDZPPN001164627} under the setting of finite element spaces. It is worth pointing out that the Lehmann--Goerisch method also requires projection-based eigenvalue bounds \cite{liu2015framework} as shown in Lemma \ref{lem:est-tau}.
Lemma \ref{thm:Lehmann-Goerisch} describes the outline of the Lehmann--Goerisch theorem used in this paper. For details, refer to \cite[\S 5.2.1]{liu2024lehmann}.

\begin{figure}[H]
\centering
\begin{tikzpicture}[scale=8]

% メッシュの節点の座標
\coordinate (n1) at (0.8000, 0.4000);
\coordinate (n2) at (0.6000, 0.3000);
\coordinate (n3) at (0.8000, 0.3000);
\coordinate (n4) at (0.8500, 0.3000);
\coordinate (n5) at (0.4000, 0.2000);
\coordinate (n6) at (0.6000, 0.2000);
\coordinate (n7) at (0.8000, 0.2000);
\coordinate (n8) at (0.8500, 0.2000);
\coordinate (n9) at (0.9000, 0.2000);
\coordinate (n10) at (0.2000, 0.1000);
\coordinate (n11) at (0.4000, 0.1000);
\coordinate (n12) at (0.6000, 0.1000);
\coordinate (n13) at (0.8000, 0.1000);
\coordinate (n14) at (0.8500, 0.1000);
\coordinate (n15) at (0.9000, 0.1000);
\coordinate (n16) at (0.9500, 0.1000);
\coordinate (n17) at (0.0000, 0.0000);
\coordinate (n18) at (0.2000, 0.0000);
\coordinate (n19) at (0.4000, 0.0000);
\coordinate (n20) at (0.6000, 0.0000);
\coordinate (n21) at (0.8000, 0.0000);
\coordinate (n22) at (0.8500, 0.0000);
\coordinate (n23) at (0.9000, 0.0000);
\coordinate (n24) at (0.9500, 0.0000);
\coordinate (n25) at (1.0000, 0.0000);

% メッシュの要素の描画
\draw (n1) -- (n2) -- (n3) -- cycle;
\draw (n2) -- (n5) -- (n6) -- cycle;
\draw (n3) -- (n2) -- (n6) -- cycle;
\draw (n3) -- (n6) -- (n7) -- cycle;
\draw (n5) -- (n10) -- (n11) -- cycle;
\draw (n6) -- (n5) -- (n11) -- cycle;
\draw (n6) -- (n11) -- (n12) -- cycle;
\draw (n7) -- (n6) -- (n12) -- cycle;
\draw (n7) -- (n12) -- (n13) -- cycle;
\draw (n10) -- (n17) -- (n18) -- cycle;
\draw (n11) -- (n10) -- (n18) -- cycle;
\draw (n11) -- (n18) -- (n19) -- cycle;
\draw (n12) -- (n11) -- (n19) -- cycle;
\draw (n12) -- (n19) -- (n20) -- cycle;
\draw (n13) -- (n12) -- (n20) -- cycle;
\draw (n13) -- (n20) -- (n21) -- cycle;
\draw (n1) -- (n3) -- (n4) -- cycle;
\draw (n3) -- (n7) -- (n8) -- cycle;
\draw (n3) -- (n8) -- (n4) -- cycle;
\draw (n4) -- (n8) -- (n9) -- cycle;
\draw (n7) -- (n13) -- (n14) -- cycle;
\draw (n7) -- (n14) -- (n8) -- cycle;
\draw (n8) -- (n14) -- (n15) -- cycle;
\draw (n8) -- (n15) -- (n9) -- cycle;
\draw (n9) -- (n15) -- (n16) -- cycle;
\draw (n13) -- (n21) -- (n22) -- cycle;
\draw (n13) -- (n22) -- (n14) -- cycle;
\draw (n14) -- (n22) -- (n23) -- cycle;
\draw (n14) -- (n23) -- (n15) -- cycle;
\draw (n15) -- (n23) -- (n24) -- cycle;
\draw (n15) -- (n24) -- (n16) -- cycle;
\draw (n16) -- (n24) -- (n25) -- cycle;

% メッシュの節点の描画
\fill (n1) circle (0.25pt);
\fill (n2) circle (0.25pt);
\fill (n3) circle (0.25pt);
\fill (n4) circle (0.25pt);
\fill (n5) circle (0.25pt);
\fill (n6) circle (0.25pt);
\fill (n7) circle (0.25pt);
\fill (n8) circle (0.25pt);
\fill (n9) circle (0.25pt);
\fill (n10) circle (0.25pt);
\fill (n11) circle (0.25pt);
\fill (n12) circle (0.25pt);
\fill (n13) circle (0.25pt);
\fill (n14) circle (0.25pt);
\fill (n15) circle (0.25pt);
\fill (n16) circle (0.25pt);
\fill (n17) circle (0.25pt);
\fill (n18) circle (0.25pt);
\fill (n19) circle (0.25pt);
\fill (n20) circle (0.25pt);
\fill (n21) circle (0.25pt);
\fill (n22) circle (0.25pt);
\fill (n23) circle (0.25pt);
\fill (n24) circle (0.25pt);
\fill (n25) circle (0.25pt);

% (n10) -- (n11) の長さが h であることを示す丸い線
\draw[black, decoration={brace, amplitude=5pt, raise=2pt, mirror}, decorate] (n17) -- node[below=4pt] {$h$} (n18);

% 三角形の頂点の座標を併記
\node[left] at (n17) {$(0, 0)$};
\node[right] at (n25) {$(1, 0)$};
\node[above right] at (n1) {$p$};
\end{tikzpicture}
\caption{Right triangular mesh}
\label{fig:right-triangular-mesh}
\end{figure}

\begin{lemma}\label{thm:Lehmann-Goerisch}
% \begin{itemize}
Given approximate eigenfunctions $\hat{u}_i \in V_h^{\CG}$ ($i=1,\cdots,n)$, let $\hat w_i \in H(\mbox{div}, T)$ be one of the functions satisfying
\begin{equation}
\label{eq:w_term_in_LH_method}
(\hat w_i, \nabla v) = (\hat u_i, v)\quad \forall v \in H^1_0(T).    
\end{equation}
Let $\lambda_{n,h}$ be the eigenvalue obtained by solving the variational equation in $V_h^{\CG}$. 
Let $\rho $ be a lower bound of the $(n+1)$th eigenvalue, i.e., $\rho \leq \lambda_{N+1}$, under the condition $ \rho> \lambda_{n,h}$.
    % \item[(A1)] Let $\hat w_i\in \mbox{RT}^p_h(T) ~(i=1,\cdots N)$ be the solution to the following variational problem:
%
 %   
    % \item[(A2)] Suppose $\rho > 0$ is a lower bound of the $(N+1)$th eigenvalue, i.e., $\rho \leq \lambda_{N+1}$.
%
    % \item[(A3)] 
    Define $n \times n$ matrices $A_0, A_1, A_2, A, B$ by
    \begin{equation}
    A_0 := \left( (\nabla\hat{u}_i,\nabla\hat{u}_j)_{T^0} \right)_{i,j=1}^{n}, \quad A_1 := \left( (\hat{u}_i, \hat{u}_j)_{T^0} \right)_{i,j=1}^{n}, \quad A_2 := \left( (\hat w_i, \hat w_j)_{T^0} \right)_{i,j=1}^{n},
    \end{equation}
    \begin{equation}
    A := A_0 - \rho A_1, \quad B := A_0 - 2\rho A_1 + \rho^2 A_2.
    \end{equation}
 Suppose $B$ is positive definite, which can be numerically validated directly.
    %
    % \item[(A5)] 
    Let $\tau_1 \leq \tau_2 \leq \cdots \leq \tau_n$ be the eigenvalues of the generalized eigenvalue problem $Az = \mu Bz$. 
    %Let $q$ be the number of negative eigenvalues.
% \end{itemize}
Under the above settings, the following lower eigenvalue bounds hold:
\begin{equation}
\lambda_{k} \geq \rho - \frac{\rho}{1 - \tau_{N+1-k}} \quad (1 \leq k \leq n).
\end{equation}
\end{lemma}

An easy-to-compute selection of $\widehat{w}_i$ is to seek $\widehat{w}_i$ in the Raviart--Thomas FEM space $\mbox{RT}_h$, which provides an approximation to the $H(\mbox{div})$ space. The Raviart--Thomas space $\mbox{RT}_h^d$ of degree $d$ is defined as follows:
\begin{equation}
    \mbox{RT}_h^d :=\{p_h\in H(\mbox{div})~|~p_h = (a_K, b_K)+c_k \mathbf{x}, a_K, b_K, c_K \in P^d(K), \mbox{ for each } K \in \mathcal{T}^h \}.
\end{equation}
The following function space $X_h^d$ of discontinous polynomials will be also used for determinating $\widehat{w}_i$'s.
\begin{equation}
X_h^d :=\{v_h ~|~ v_h \in P^d(K) \mbox{ for each } K \in \mathcal{T}^h \}.
\end{equation}

Note that the solution to $\eqref{eq:w_term_in_LH_method}$ is not unique.
An optimal selection of $\widehat{w}_i$ can be determined by solving the problem: Find $(\hat p_i,\hat g_i)\in\mbox{RT}^d_h\times X_h^{d}$ such that
    \begin{equation}
        \begin{cases}
            (\hat p_i,\hat q)+(\hat g_i,\mbox{div } \hat q)=0~~\forall \hat q\in \mbox{RT}^d_h,\\
            (\mbox{div }\hat p_i, \hat f) + (\hat u_i,\hat f)=0 ~~\forall \hat f\in X_h^{d}.
        \end{cases}
    \end{equation}
In the computation of high-precision gienvalue  bounds, the degree $d$ is taken as $d=5$.

        \medskip\medskip\medskip
             
        Next, we recall the estimations concerning the perturbation of eigenvalues, which is needed to compute the eigenvalue bounds over triangular domains with a given range of parameters.
								
	\medskip
								
	\textbf{Perturbation of functions with respect to variation of triangles}

        \medskip
 
	For $p=(x,y)\in \Omega$ and its perturbation $\tilde p=(\tilde x,\tilde y)$, consider the linear transformation $\Phi:T^{p}\to T^{\tilde p}$. The representation matrix $S$ for $\Phi$ is given by
	\begin{equation}
		S:=
		\begin{pmatrix}
			1 & (\tilde x-x)/y \\
			0 & \tilde y/y     
		\end{pmatrix}.
	\end{equation}
								
	For $u\in H^1_0(T^p)$, define $\tilde{u}$ over $H^1_0(T^{\tilde p})$ by
	$\tilde{u}=u\circ\Phi^{-1}$.
	Denote by $S^\intercal$ the transpose of $S$.
	It holds that $\nabla\tilde{u}(\tilde x,\tilde y)=S^{-\intercal}
	\nabla u(x,y)$.
	Let  
	$\lambda_{\min}(\cdot) $ and $\lambda_{\max}(\cdot)$ denote the minimum and the maximum eigenvalues of a given square symmetric matrix, respectively.
								
	\medskip
								
	Let us recall the result about the perturbation of eigenvalues with respect to perturbation of eigenvalues, as previously well investigated in \cite{liu2015framework} and formulated in the general case in \cite{endo2023shape}.

        \begin{lemma}[Lemma Appendix A.2 of \cite{endo2023shape}]
        \label{lem:eig-perturbation}
            For the $i$-th $(i=1,2,\cdots)$ eigenvalue at $p=(x,y)$, we have
            \begin{equation}\label{eq:perturbation}
                \lambda_{\min}\left(S^{-1}S^{-\intercal}\right)\cdot\lambda^{p}_i\leq \lambda^{\tilde p}_i\leq\lambda_{\max}\left(S^{-1}S^{-\intercal}\right)\cdot\lambda^{p}_i.
            \end{equation}
        \end{lemma}

        The coefficients in the above perturbation estimation formulas for perturbations of $p$ in the $x$-direction and $y$-direction are described in the following remark.
        \begin{remark}\label{eq:perturbation-monotonicity}
            In case $y=\tilde y$, we have
            \begin{gather}\label{eq:perturbation_special_1}
            \lambda_{\min}\left(S^{-1}S^{-\intercal}\right)=\frac{1}{2}\left(2+\frac{(x-\tilde x)^2}{y^2}\right)-\frac{|x-\tilde x|}{2y}\sqrt{4+\frac{(x-\tilde x)^2}{y^2}}~(=:m_x(p,\tilde p)),\\
            \lambda_{\max}\left(S^{-1}S^{-\intercal}\right)=\frac{1}{2}\left(2+\frac{(x-\tilde x)^2}{y^2}\right)+\frac{|x-\tilde x|}{2y}\sqrt{4+\frac{(x-\tilde x)^2}{y^2}}~(=:M_x(p,\tilde p)).
            \end{gather}
            Note that, $m_x(p,\tilde p)$ (resp. $M_x(p,\tilde p)$) is monotonically decreasing (resp. increasing) with respect to $|x-\tilde x|$.
            
            Also, if $x=\tilde x$ and $\tilde y\leq y$, then we have
            \begin{gather}\label{eq:perturbation_special_2}
            \lambda_{\min}\left(S^{-1}S^{-\intercal}\right)=1~(=:m_y(p,\tilde p)),~~
            \lambda_{\max}\left(S^{-1}S^{-\intercal}\right)=\frac{y^2}{\tilde y^2}~(=:M_y(p,\tilde p)).
            \end{gather}
        \end{remark}        
							
	\subsection{Step 1: Estimation of the difference quotients over nearly equilateral triangles}
\begin{figure}[H]
  \centering
  \includegraphics[keepaspectratio, scale=6]{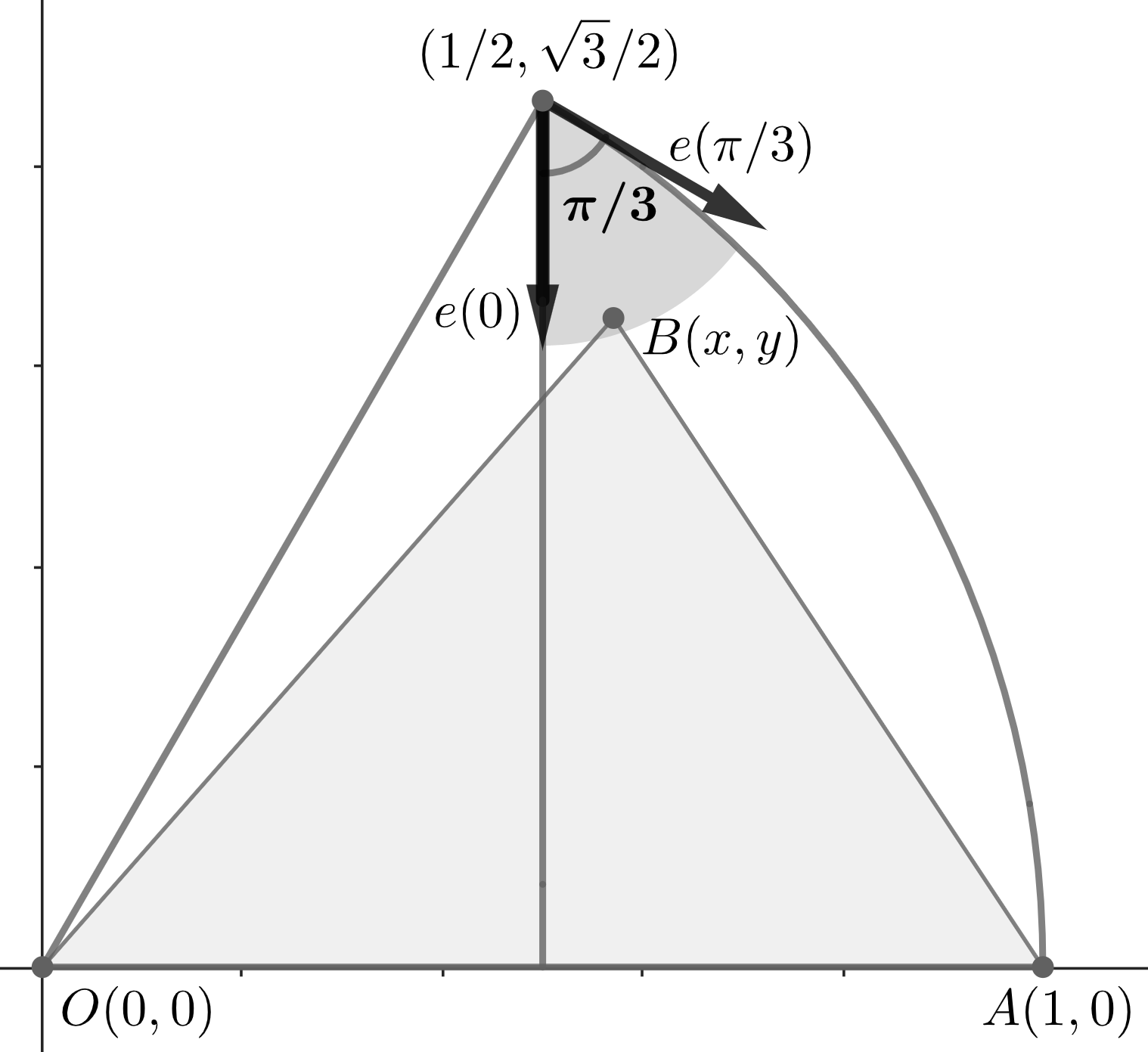}
  \caption{The case $p\in\Omega_{up}$\label{fig:turn-around}}
\end{figure}

Let $T^p$ be a perturbation of the equilateral triangle $T^{p_0}$. In this section, we explicitly estimate the range of the difference quotients of the second and third eigenvalues for triangles slightly perturbated from the equilateral triangle. Consider the closed ball $B(p_0;\varepsilon)$ with $\varepsilon=10^{-5}$. We estimate the difference quotient of the eigenvalues, $D_t\lambda_i(p_0,p)$, for $p\in B(p_0;\varepsilon)$ that contains $\Omega_{up}$. For $\delta\in[0,\pi/3]$, define the direction of derivatives as $e(\delta):=(\sin\delta ,-\cos\delta )$, as illustrated in Figure \ref{fig:turn-around}.

Recall that the Dirichlet eigenvalues over the equilateral triangle have the following explicit values \cite{mccartin2003eigenstructure}:
\begin{equation}\label{eq:explicit-values-of-eigenvalues}
  \lambda_{1}^{p_0}=\frac{16}{3}\pi^2,~\lambda_{2}^{p_0}=\lambda_{3}^{p_0}=\frac{112}{9}\pi^2,~\lambda_{4}^{p_0}=\frac{64}{3}\pi^2~~\mbox{for }p_0=(1/2,\sqrt{3}/2).
\end{equation}

\textbf{Estimation of the difference quotients}

\begin{figure}[H]
    \centering
    \includegraphics[keepaspectratio, scale=0.6]{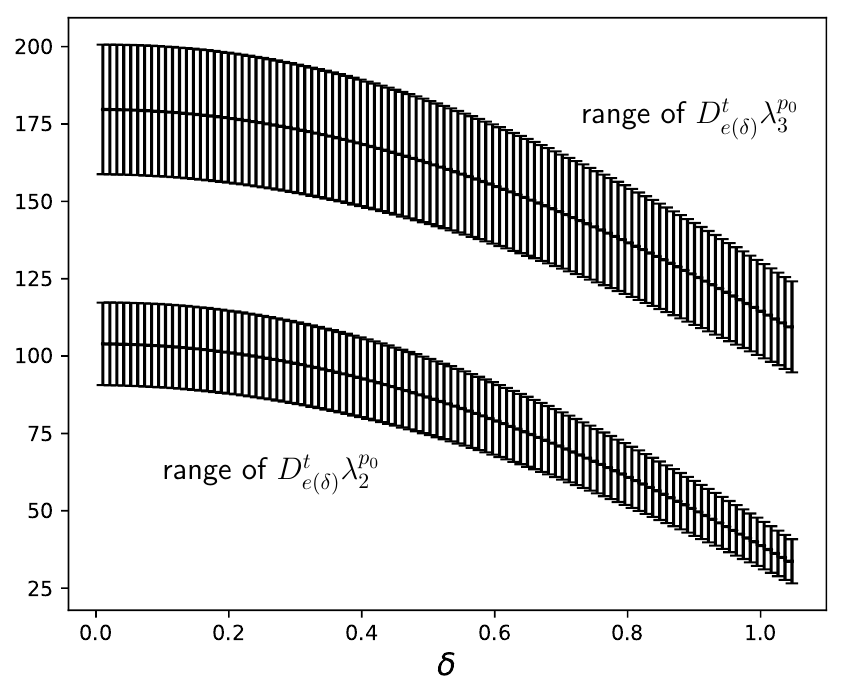}
    \caption{\label{fig:ranges-of-derivatives} The obtained ranges of $D_{e(\delta)}^t\lambda_2^{p_0}$ and $D_{e(\delta)}^t\lambda_3^{p_0}$}
\end{figure}

Let us estimate the range of the difference quotient $D_t\lambda_i(p_0,p)$ for $p\in B(p_0;\varepsilon)\backslash\{p_0\}$ using Algorithm \ref{algorithm-1}. For each $\lambda_i$ $(i=2,3)$, define $$D_{e(\delta)}^t\lambda_i^{p_0}:=D_t\lambda_i(p_0,p_0+te(\delta)),\quad  ~t\in(0,\varepsilon], ~ \delta\in[0,\pi/3].
$$
Divide the interval $[0,\pi/3]$ into subintervals $[\delta_k,\delta_{k+1}]$, where $\delta_k:=k{\pi}/3000$ for $k=0,\ldots,999$. Figure \ref{fig:ranges-of-derivatives} shows the estimated ranges of $D_{e(\delta)}^t\lambda_2^{p_0}$ and $D_{e(\delta)}^t\lambda_3^{p_0}$, and Tables \ref{table:problem-1-quantities-1} and \ref{table:problem-1-quantities-2} provide the values of related quantities.
										
		As a result of the computation, it follows that for all $t\in (0,\varepsilon]$ and $\delta\in [0,\pi/3]$, we have
$$
  D_{e(\delta)}^t\lambda_2^{p_0}<D_{e(\delta)}^t\lambda_3^{p_0}, \mbox{ i.e.},  \frac{\lambda_2^{p} - \lambda_2^{p_0}}{t} <  \frac{\lambda_3^{p} - \lambda_3^{p_0}}{t}, ~\mbox{implying }
  \lambda_2^p<\lambda_3^p\quad 
 \forall p\in\Omega_{up}\backslash\{p_0\}.
$$

\begin{algorithm}[H]
\caption{Rigorous bounds for difference quotients $D_{e(\delta)}^t\lambda_i^{p_0}$ ($i=2,3$) over subinterval $[\delta_k,\delta_{k+1}]$ of $[0,\pi/3]$}
\label{algorithm-1}
\begin{algorithmic}[1]
\Require Intervals for perturbation magnitude $(0,\varepsilon]$ and direction $[\delta_k,\delta_{k+1}]$.
\Ensure Intervals  $I_{k,i} (i=2,3)$ as enclosures of  difference quotients $D_{e(\delta)}^t\lambda_i^{p_0}$ ($i=2,3$), respectively,  for all $t\in(0,\varepsilon]$ and $\delta\in[\delta_k,\delta_{k+1}]$. That is, 
\begin{equation*}
    \{D_{e(\delta)}^t\lambda_i^{p_0}~|~t \in (0,\varepsilon],~\delta \in [\delta_k,\delta_{k+1}]\}\subset I_{k,i}~~(i=2,3).
\end{equation*}

\medskip

\State \textbf{(Preparation)} Estimate the range of $\lambda_i^{p}$ ($i=2,3,4$) over the ball $B(p_0;\varepsilon)$ using the explicit values from  \eqref{eq:explicit-values-of-eigenvalues} and the perturbation results from Lemma \ref{lem:eig-perturbation}.
\begin{equation}\label{eq:equi-pert-est}
    \lambda_{\min}\left([S^{-1}S^{-\intercal}]\right)\cdot\lambda^{p_0}_i\leq \lambda^{p}_i\leq\lambda_{\max}\left([S^{-1}S^{-\intercal}]\right)\cdot\lambda^{p_0}_i,
\end{equation}
where the interval matrix $[S^{-1}S^{-\intercal}]$ is defined by
\begin{align*}
     [S^{-1}S^{-\intercal}]&:=\{S^{-1}S^{-\intercal}~|~S:T^{p_0}\to T^{p} \mbox{ being a linear map},~~p\in B(p_0;\varepsilon)\}.
\end{align*}
Here, $\lambda_{\min}(\cdot)$ and $\lambda_{\max}(\cdot)$ take the minimum and the maximum eigenvalues of an interval matrix, respectively.

\State \textbf{(FEM Basis)} Calculate approximate eigenfunctions $\hat u_2, \hat u_3 \in H^1_0(T^{p_0})$ for the equilateral triangle using FEM basis functions $\{\psi_k\}_{k=1}^m$. Let $\mathbf{c}_2, \mathbf{c}_3 \in \mathbb{R}^m$ be the coefficient vectors of $\hat u_2, \hat u_3$.  

\State \textbf{(Construct Interval Matrices)} 
Let us compute the matrices in the difference quotient formula \eqref{def:approx-matrices}.
Define $2\times 2$ interval matrix $[P]$ by
\begin{equation*}
    [P]_{ij}:=\{(P_t^{e(\delta)})_{ij}~|~t \in (0,\varepsilon],~\delta \in [\delta_k,\delta_{k+1}]\}~~(i,j=1,2).
\end{equation*}
Calculate interval matrices $[A]$ and $B$ using interval arithmetic:
\begin{equation*}
    [A]_{ij}:=\int_{T^{p_0}}[P]\nabla\psi_i\cdot\nabla\psi_j~dx,~~[B]_{ij}:=\int_{T^{p_0}}\psi_i\psi_j~dx~~(i,j=1,\cdots,m).
\end{equation*}
The elements of  $[\widehat M_t]$ and $[\widehat N_t]$ are created by
\begin{equation*}
    (\widehat M_t)_{ij} = \mathbf{c}_{i+1}^T [A] \mathbf{c}_{j+1},~~(\widehat N_t)_{ij} = \mathbf{c}_{i+1}^T [B] \mathbf{c}_{j+1}~~(i,j=1,2).
\end{equation*}

\end{algorithmic}
\end{algorithm}

\begin{algorithm}[H]
\ContinuedFloat
\caption{Rigorous bounds for difference quotients $D_{e(\delta)}^t\lambda_i^{p_0}$ ($i=2,3$) over $[\delta_k,\delta_{k+1}]$ (continued)}
\begin{algorithmic}[1]
\setcounter{ALG@line}{3}

\State \textbf{(Estimate Eigenspace Errors)} 
Note that $ \widehat{E}
=\mbox{span}\{\hat u_2, \hat u_3\},  {E}
=\mbox{span}\{u_2, u_3\} $, ${\widetilde{E}_t}
=\mbox{span}\{\tilde{u}^t_2, \tilde{u}^t_3\} $.
Compute uniform upper bounds for the error terms $\mathrm{Err}_F$ and $\mathrm{Err}_b$ from Lemma \ref{lem:M-star-N-star-estimation}. This requires rigorously bounding the directed distances $\overline{\delta}_b(\widehat E,E)$ and $\overline{\delta}_b(E,\widetilde E_t)$ by applying the framework of Lemma \ref{lem:eigenvec-estimation-algorithm-1}.
\begin{itemize}
    \item[\textbf{a)}] To bound the \textbf{FEM approximation error} $\overline{\delta}_b(\widehat E,E)$,    
    apply Lemma \ref{lem:eigenvec-estimation-algorithm-1} with the following settings:
    $$
    E_1 = \mathrm{span}\{u_1\}, \quad 
    E_2 =  \mathrm{span}\{u_2, u_3\},\quad E^h_1 = \mathrm{span}\{\hat{u}_1\},\quad
    E^h_2 = \mbox{span}\{\hat u_2, \hat u_3\}    $$
    along with explicit value of eigenvalue $\lambda^{p_0}_i$ ($i=1,\cdots, 4)$ given in \eqref{eq:explicit-values-of-eigenvalues} and $\rho:=\lambda^{p_0}_4$.

    \item[\textbf{b)}] To bound the \textbf{domain perturbation error} $\overline{\delta}_b(E,\widetilde E_t)$,    
    apply Lemma \ref{lem:eigenvec-estimation-algorithm-1} with the following settings:
    $$
    E_1 = \mathrm{span}\{u_1\}, \quad 
    E_2 =  \mathrm{span}\{u_2, u_3\},\quad E^h_1 = \mathrm{span}\{\tilde u_1^t\},\quad
    E^h_2 = \mbox{span}\{\tilde u_2^t, \tilde u_2^t\}    $$
    along with explicit value of eigenvalue $\lambda^{p_0}_i$ ($i=1,\cdots, 4)$ given in \eqref{eq:explicit-values-of-eigenvalues} and $\rho:=\lambda^{p_0}_4$, and
    $$\lambda_{3,h} := \max_{\tilde u\in \widetilde E_t}\frac{\|\nabla\tilde u\|_{T_0}^2}{\|\tilde u\|_{T_0}^2}.$$
    
    \item[\textbf{c)}] \textbf{Final Error Calculation:}
    The computed bounds for $\delta_b$ are converted to $\overline{\delta}_b$ using the relation in Lemma \ref{lem:bbar-b-relation}. These are then used to evaluate the intermediate error quantities:
    \begin{align*}
        \hat\eta^2 &:= \lambda_N\overline{\delta}_b^2(\widehat E,E)+\hat{\lambda}_N-\lambda_n,~~\\
        \tilde\eta^2 &:= \lambda_N\overline{\delta}_b^2(E,\widetilde E_t)+\tilde{\lambda}_N^t-\lambda_n,
    \end{align*}
    Here, $\tilde\lambda_N^t$ is bounded by $\tilde\lambda_N^t \le \|S_tS_t^{\intercal}\|^2\lambda_N^t$.
    These values are then substituted into the formulas from Lemma \ref{lem:M-star-N-star-estimation} to obtain the final error bounds $\mathrm{Err}_F$ and $\mathrm{Err}_b$.
\end{itemize}

\State \textbf{(Construct Guaranteed Interval Matrices)} Inflate the approximate matrices with the computed error bounds to form guaranteed interval matrices $[M_t]$ and $[N_t]$:
\begin{align*}
    [M_t] &= [\widehat M_t] + [-\mathrm{Err}_F, \mathrm{Err}_F] \\
    [N_t] &= [\widehat N_t] + [-\mathrm{Err}_b, \mathrm{Err}_b]
\end{align*}
These matrices are guaranteed to contain the true matrices $M_t^*$ and $N_t^*$.

\State \textbf{(Solve Interval Eigenvalue Problem)} Solve the generalized interval eigenvalue problem $[M_t]\sigma=\mu[N_t]\sigma$. This is performed using a verified numerical method, such as those implemented in the INTLAB library \cite{Rump1999intlab}, to obtain the rigorous interval enclosures $I_{k,i}$ for the true difference quotients.
\end{algorithmic}
\end{algorithm}

\begin{remark}
Note that the specific value $\varepsilon=10^{-5}$ is determined based on the requirements of verified computation. If the value of $\varepsilon$ is set too large, the value of $\overline{\delta}_b(E,\widetilde E_t)$ increases, along with the values of $\mbox{Err}_F(P^e_t;E,\widehat E)$ and $\mbox{Err}_b(P^e_t;E,\widehat E)$, causing the intervals $I_{k,i}(k=2,3)$, which evaluate $\mu_k(k=2,3)$, to overlap. Conversely, if the value of $\varepsilon$ is set too small, $\lambda_2^p$ and $\lambda_3^p$ become very close in the neighborhood of $p_0$, making it difficult to separate them in Step 2. We chose the value $\varepsilon=10^{-5}$ so that $I_{k,i}~(k=2,3)$ do not overlap and eigenvalue separation is possible in $\Omega_{down}$.
\end{remark}

\begin{figure}[H]
    \centering
    \begin{minipage}[t]{0.48\textwidth}
        \centering
        \captionof{table}{\label{table:problem-1-quantities-1} The obtained range of $D_{e(0)}^t\lambda_2^{p_0}$, $D_{e(0)}^t\lambda_3^{p_0}$ and related quantities}
        \begin{tabular}{cl}
            \hline
            \rule[-2mm]{0mm}{6mm}{}
            $I_{2,0}$ & $[89.5557, 118.3308]$ \\
            $I_{3,0}$ & $[162.9557, 196.4519]$ \\
            \hline
            \rule[-2mm]{0mm}{6mm}{}
            $\hat\mu_2$                         & $\approx 103.8982$       \\
            \rule[-2mm]{0mm}{6mm}{}
            $\hat\mu_3$                         & $\approx 179.6588$       \\
            \rule[-2mm]{0mm}{6mm}{}
            $\mbox{Err}_F(P^e_t;E,\widehat E)$  & $\leq 2.7760$           \\
            \rule[-2mm]{0mm}{6mm}{}
            $\mbox{Err}_b(P^e_t;E,\widehat E)$  & $\leq 0.0078$            \\
            \hline
        \end{tabular}
    \end{minipage}
    \hspace{0.04\textwidth} % Adjust space between minipages as needed
    \begin{minipage}[t]{0.48\textwidth}
        \centering
        \captionof{table}{\label{table:problem-1-quantities-2} The obtained range of $D_{e(\frac{\pi}{3})}^t\lambda_2^{p_0}$, $D_{e(\frac{\pi}{3})}^t\lambda_3^{p_0}$ and related quantities}
        \begin{tabular}{cl}
            \hline
            \rule[-2mm]{0mm}{6mm}{}
            $I_{2,999}$ & $[21.1589,  45.0331]$  \\
            $I_{3,999}$ & $[93.9566, 123.7559]$ \\
            \hline
            \rule[-2mm]{0mm}{6mm}{}
            $\hat\mu_2$                         & $\approx 33.0085$        \\
            \rule[-2mm]{0mm}{6mm}{}
            $\hat\mu_3$                         & $\approx 108.7637$      \\
            \rule[-2mm]{0mm}{6mm}{}
            $\mbox{Err}_F(P^e_t;E,\widehat E)$  & $\leq 2.6145$           \\
            \rule[-2mm]{0mm}{6mm}{}
            $\mbox{Err}_b(P^e_t;E,\widehat E)$  & $\leq 0.0098$            \\
            \hline
        \end{tabular}
    \end{minipage}
\end{figure}

             % [\LIU{To be removed.}]
						      	% Since we have $D_{e(\delta)}^t\lambda_2^{p_0}<D_{e(\delta)}^t\lambda_3^{p_0}$, it follows that $\lambda_2^{p}<\lambda_3^{p}$ over $\Omega_{up}\backslash\{0\}$. ($\star$)
						      
\subsection{Step 2: Eigenvalue bounds by finite element method}

We present an algorithm for the verified numerical computations required in Step 2, along with the computational results.

In Step 2, we cover $\ENDO{\Omega_{down}^{(1)}}$ and $\ENDO{\Omega_{down}^{(2)}}$ with small rectangles $R_{ij}$, and we compute eigenvalue bounds at each node of the rectangles using direct eigenvalue estimation (Lemmas \ref{lem:est-tau} and \ref{thm:Lehmann-Goerisch}). For $\ENDO{\Omega_{down}^{(1)}}$, in which the height of the triangle $T^p$ is greater than $\tan(\pi/14)$, we apply the perturbation estimates (Lemma \ref{lem:eig-perturbation}) to obtain a uniform bound on each rectangle $R_{ij}$. In the case of $\ENDO{\Omega_{down}^{(2)}}$, where the triangle $T^p$ is nearly degenerate, we employ the domain monotonicity property of Dirichlet eigenvalues to obtain a uniform bound on each rectangle $R_{ij}$. Details are provided below.

\paragraph{Case of $\ENDO{\Omega_{down}^{(1)}}$:}

To compute the eigenvalue bounds \( \lambda_k \) over $\ENDO{\Omega_{down}^{(1)}}$ for \( k = 2,3 \), we partition the set \(\ENDO{\Omega_{down}^{(1)}} \) into small intervals. For indices \( i = 1,2,\cdots,5000 \) and \( j = 1,2,\cdots,1020 \), define nodes $\{(x_{i},y_{j})\}$ by

\begin{gather*}
  (x_i,y_j)=
  \begin{cases}
     ( 0.5+10^{-6}(i-1), ~\frac{\sqrt{3}}{2}-10^{-5}-10^{-6}(j-1) ) & (1\leq j\leq 91)\\
     ( 0.5+10^{-5}(i-1),  ~ y_{91}-10^{-5}(j-91)) & (92\leq j\leq 181)\\
     ( 0.5+10^{-4}(i-1), ~ y_{181}-10^{-4}(j-181)) & (182\leq j\leq 271)\\
     ( 0.5+10^{-3}(i-1), ~ y_{271}-10^{-3}(j-271)) & (272\leq j\leq 1020)
  \end{cases}.
\end{gather*}

Let $R_{ij}$ be the rectangle with vertices $(x_{i},y_{j}),(x_{i},y_{j+1}),(x_{i+1},y_{j+1}),(x_{i+1},y_{j})$, and let \( \Omega_{ij}^{(1)} := R_{ij} \cap \Omega \).
Note that, from the estimation \eqref{eq:perturbation} and the monotonicity described in Remark \ref{eq:perturbation-monotonicity}, we have
\begin{equation}\label{eq:perturbation-est-in-R}
    m_x((x_{i+1},y_{j}),(x_{i},y_{j}))\lambda_k^{(x_{i+1},y_j)}\leq\lambda_k^p\leq
    M_x((x_{i},y_{j+1}),(x_{i+1},y_{j+1}))\lambda_k^{(x_i,y_{j+1})}
\end{equation}
for $p\in R_{ij},~k=2,3$.

To compute the uniform bounds of the eigenvalues $\lambda_k$ (for $k=2,3$) in each $\Omega_{ij}^{(1)}$, we employ the following algorithm:
                    
\begin{algorithm}[H]
    \caption{\label{algorithm-2} Upper and lower bounds of $\lambda_k^p$ over the domain $\Omega$.}
    \KwData{Subdivision $\Omega_{ij}^{(1)}$}
    \KwResult{$I_{ij}$ as an interval that uniformly encloses $\lambda_k^p$ ($k=2,3$) for all $p=(x,y) (\in \Omega_{ij}^{(1)})$.}
    \textbf{Procedure:} 
    \begin{enumerate}
        \item [1.] Obtain an upper bound of $\lambda_k^{(x_i,y_{j+1})}$ as $\overline\lambda_k^{(x_i,y_{j+1})}$, and a lower bound of $\lambda_k^{(x_{i+1},y_j)}$ as $\underline\lambda_k^{(x_{i+1},y_j)}$ by Lemma \ref{lem:est-tau} and \ref{thm:Lehmann-Goerisch}.
        \item [2.] Obtain $I_{ij}=[\underline{I}_{ij},\overline{I}_{ij}]$ as uniform bounds of $\lambda_k^p$ over $R_{ij}$  by letting 
        \begin{align}
            \underline{I}_{ij}&=m_x((x_{i},y_{j}),(x_{i+1},y_{j}))\underline\lambda_k^{(x_{i+1},y_j)},\\
            \overline{I}_{ij}&=M_x((x_{i+1},y_{j+1}),(x_{i},y_{j+1}))\overline\lambda_k^{(x_i,y_{j+1})}
        \end{align}
        using the estimate \eqref{eq:perturbation-est-in-R}.
    \end{enumerate}
\end{algorithm}

The obtained lower bound of the gap $\lambda_{3}^p-\lambda_2^p$ over $\ENDO{\Omega_{down}^{(1)}}$ is as follows:
\begin{equation}
    \min_{p\in \ENDO{\Omega_{down}^{(1)}}}(\lambda_{3}^p-\lambda_2^p)\geq 8.4814\cdot 10^{-4}
\end{equation}

\begin{remark}
Let us confirm the indispensability of the Lehmann--Goerisch method in separating eigenvalues.
When $T^p~(p\in\ENDO{\Omega_{down}^{(1)}})$ is close to an equilateral triangle, that is, when $\|p-p_0\|_2\approx 10^{-5}$, it is not possible to separate $\lambda_2^p$ and $\lambda_3^p$ using only the estimation from Lemma \ref{lem:est-tau} without employing the Lehmann-Goerisch Method.
For example, for $p=(1/2,\sqrt{3}/2-10^{-5})$, the estimation obtained using the Lehmann-Goerisch Method is
\begin{equation}
122.8225\leq\lambda_2\leq 122.8231,~~122.8232\leq\lambda_3\leq 122.82393, 
\end{equation}
which successfully separates the eigenvalues, whereas the estimation using only Lemma \ref{lem:est-tau} gives
\begin{equation}
122.0788\leq\lambda_2\leq 122.9701,~~122.1173\leq\lambda_3\leq 123.0092,
\end{equation}
which results in an overlap.    
\end{remark}

\begin{remark}
For $p=(x,y)$ with $(x-0.8893)^2+(y-0.4571)^2<0.01$, we have $\lambda_3^p\approx\lambda_4^p$. In this case, to apply the Lehmann--Goerisch method, one have to use the lower bound of $\rho<\lambda_5$ and take $n=4$ in Lemma \ref{thm:Lehmann-Goerisch}. For this case, we obtain a lower bound for $\lambda_3^p$ using only Lemma \ref{lem:est-tau} rather than Lemma \ref{thm:Lehmann-Goerisch}.
\end{remark}

\paragraph{Case of $\ENDO{\Omega_{down}^{(2)}}$:}

In the case of of $p\in\ENDO{\Omega_{down}^{(2)}}$, we apply only Lemma \ref{lem:est-tau} to validate the simplicity of $\lambda_2^p$ without using Lehmann--Goerisch method, since the eigenvalues $\lambda_2^p$ and $\lambda_3^p$ are well separated.

For the parameter $p=(x,y)\in \ENDO{\Omega_{down}^{(2)}}$, we introduce the coordinate transformation
\[
r = \sqrt{x^2+y^2},\quad h = y.
\]
Due to the monotonicity of Dirichlet eigenvalues with respect to domain inclusion, the Dirichlet eigenvalue $\lambda_k^{(r,h)}$ on the triangle $T^p$ is monotonically decreasing in both $r$ and $h$; see Figure \ref{fig:monotonicity-domain}.

\begin{figure}[H]
  \begin{minipage}[c]{0.45\linewidth}
    \centering
    \small 
    \begin{tikzpicture}[scale=5.5]
  % Draw the common base line
  \draw (0,0) -- (1,0);
  
  % Draw the smaller triangle T^{p_1} (with lower r)
  % Let p_1 = (0.6,0.4)
  \draw (0.0,0) -- (1,0) -- (0.6,0.4) -- cycle;
  \node at (0.6,0.4) [above] {$p_2$};
  \node at (0.42,0.1) {$T^{p_1}$};
  
  % Draw the larger triangle T^{p_2} (with higher r) using dashed lines
  % Let p_2 = (0.4,0.266)
  \draw[thick] (0,0) -- (1,0) -- (0.4,0.266) -- cycle;
  \node at (0.4,0.266) [above] {$p_1$};
  \node at (0.6,0.266) {$T^{p_2}$};

  \node at (0,0) [below] {$(0,0)$};
  \node at (1,0) [below] {$(1,0)$};
  
\end{tikzpicture}
    \label{fig:eigenvalue_monotonicity_r}
  \end{minipage}
  ~~
  \begin{minipage}[c]{0.45\linewidth}
    \centering
    \small 
    \begin{tikzpicture}[scale=5.5]
  % Draw the common base line
  \draw (0,0) -- (1,0);
  
  % Draw the smaller triangle T^{p_1} (with lower r)
  % Let p_1 = (0.6,0.4)

  \draw[dashed] (0.6,0.0) -- (0.6,0.4) -- cycle;
  
  \draw (0.0,0) -- (1,0) -- (0.6,0.4) -- cycle;
  \node at (0.6,0.42) [right] {$p_2$};
  \node at (0.47,0.1) {$T^{p_1}$};
  
  % Draw the larger triangle T^{p_2} (with higher r) using dashed lines
  % Let p_2 = (0.4,0.266)
  \draw[thick] (0,0) -- (1,0) -- (0.6,0.266) -- cycle;
  \node at (0.65,0.26) [above] {$p_1$};
  \node at (0.5,0.26) {$T^{p_2}$};

  \node at (0,0) [below] {$(0,0)$};
  \node at (1,0) [below] {$(1,0)$};
  
\end{tikzpicture}
    \label{fig:eigenvalue_monotonicity_h}
  \end{minipage}
  \caption{Domain monotonicity with respect to $r$ (left) and $h$ (right)\label{fig:monotonicity-domain}}
\end{figure}

To cover $\ENDO{\Omega_{down}^{(2)}}$, we partition it into subsets
\[
\Omega_{ij}^{(2)} := \Big\{(x,y)\in \ENDO{\Omega_{down}^{(2)}} : \, r\in[r_i,r_{i+1}],\; h\in[h_j,h_{j+1}]\Big\}.
\]
Here, for indices \( i = 1,2,\cdots,11 \) and \( j = 1,2,\cdots,81 \), the nodes $\{(r_i, h_j)\}$ are defined as follows:
\begin{gather*}
  (r_i,h_j)=
  \begin{cases}
     ( (i-1)/20+0.5, ~\tan(\pi/60)/2+(\tan(\pi/40)-\tan(\pi/60))(j-1)/80 ) & (1\leq j\leq 41)\\
     ( (i-1)/20+0.5,   ~\tan(\pi/40)/2+(\tan(\pi/14)-\tan(\pi/40))(j-41)/80 )) & (42\leq j\leq 81)\\
  \end{cases}.
\end{gather*}
Note that the family of rectangles $\widetilde{R}_{ij}=[r_i,r_{i+1}]\times [h_j,h_{j+1}]$ forms a covering of the following set:
\begin{equation}
  \left\{\, (r,h) \in \mathbb{R}^2 \;\middle|\;
  \begin{aligned}
    & (x, y) \in \Omega_{\mathrm{down}}^{(2)},~~ r = \sqrt{x^2 + y^2}, ~~ h = y
  \end{aligned}
  \right\}
\end{equation}
Using the monotonicity of $\lambda_k(r,h)$ in $r$ and $h$, for any $(r,h)\in \widetilde{R}_{ij}$, we obtain
\begin{equation}\label{eq:monotonicity-estimate}
\lambda_k^{(r_{i+1},h_{j+1})} \leq \lambda_k^{(r,h)} \leq \lambda_k^{(r_i,h_j)}
\end{equation}
for $k=2,3$.

Based on \eqref{eq:monotonicity-estimate}, we use the following algorithm to compute uniform eigenvalue bounds over each $\Omega_{ij}^{(2)}$:

\begin{algorithm}[H]
    \caption{\label{algorithm-4} Upper and Lower Bounds of $\lambda_k^p$ over $\ENDO{\Omega_{down}^{(2)}}$.}
    \KwData{A covering $\{\Omega_{ij}^{(2)}\}$ of $\ENDO{\Omega_{down}^{(2)}}$ as defined above.}
    \KwResult{An interval $I_{ij}=[\underline{I}_{ij},\overline{I}_{ij}]$ that uniformly encloses $\lambda_k^p$ for all $p\in \Omega_{ij}^{(2)}$, for $k=2,3$.}
    \textbf{Procedure:}
    \begin{enumerate}
        \item [1.] Compute an upper bound for $\lambda_k^{(r_i,h_j)}$, denoted by $\overline{\lambda}_k^{(r_i,h_j)}$, and a lower bound for $\lambda_k^{(r_{i+1},h_{j+1})}$, denoted by $\underline{\lambda}_k^{(r_{i+1},h_{j+1})}$, using Lemmas \ref{lem:est-tau} and \ref{thm:Lehmann-Goerisch}.
        \item [2.] Obtain $I_{ij}=[\underline{\lambda}_k^{(r_{i+1},h_{j+1})},\overline{\lambda}_k^{(r_i,h_j)}]$ as uniform bounds of $\lambda_k^p$ over $\ENDO{\Omega_{down}^{(2)}}$.
    \end{enumerate}
\end{algorithm}

The computed lower bound for the gap $\lambda_3^p-\lambda_2^p$ over $\ENDO{\Omega_{down}^{(2)}}$ is given by
\[
\min_{p\in \ENDO{\Omega_{down}^{(2)}}}\Big(\lambda_3^p-\lambda_2^p\Big)\geq 27.87.
\]
			
\section{Conclusions}
We have established a difference quotient formula for eigenvalues and derived an error estimate for its numerical approximation. Using these results, we have provided a computer-assisted proof for the simplicity of the second Dirichlet eigenvalue for triangles with minimum normalized height greater than $\tan(\pi/60)/2$.  In conjunction with the results from our previous work~\cite{endo2025second1} for triangles with  minimum normalized height less than or equal to $\tan(\pi/60)/2$, we arrive at Theorem \ref{thm:main-conjecture}.

\section*{Acknowledgement}
Both authors are supported by Japan Society for the Promotion of Science. The first author is supported by JSPS KAKENHI Grant Number JP24KJ1170. The last author is supported by JSPS KAKENHI Grant Numbers JP20KK0306, JP22H00512, JP24K00538 and JP21H00998.

                                % https://www.jsps.go.jp/j-grantsinaid/16_rule/syaji.html
                                
             %                    The last author is supported by Fund for the Promotion of Joint International Research (Fostering Joint International Research (A)) 20KK0306, Grant-in-Aid for Scientific Research (A) 22H00512, 
						      	% Grant-in-Aid for Scientific Research (B) 20H01820, 21H00998, and Grant-in-Aid for Scientific Research (C) 21K03367. 													
						      	\bibliographystyle{elsarticle-num} 
						      	\bibliography{references}
						      						 
\end{document}